\documentclass[12pt,twoside]{article}
\usepackage{times}
\usepackage{amsmath,amssymb}
\usepackage{amsthm}
\usepackage{mathrsfs}
\usepackage{color}
\usepackage{exscale}
\usepackage{relsize}
\usepackage{color}
\usepackage{extarrows}
\usepackage{fancyhdr}
\allowdisplaybreaks

\newcommand{\pf}{\noindent {\bf Proof. \hspace{2mm}}}
\newcommand{\ef}{ \hfill $ \Box $ \vskip 3mm}
\newcommand{\be}{\begin{equation}}
\newcommand{\ee}{\end{equation}}
\newcommand{\bea}{\begin{eqnarray}}
\newcommand{\eea}{\end{eqnarray}}

\newcommand{\bR}{{\mathbb R}}
\newcommand{\bN}{{\mathbb N}}
\newcommand{\bZ}{{\mathbb Z}}

\def\nn{\nonumber}
\def\R {\Bbb R}

\def\ve{\varepsilon}
\def\la{\lambda}

\def\t{\tilde}
\def\q{\quad}
\def\qq{\qquad}
\def\th{\theta}

\def\G{\Gamma}
\def\dl{\delta}
\def\Dl{\Delta}
\def\ve{\varepsilon}

\def\lt{\left}

\def\rt{\right}

\def\i{\infty}

\def\supp{\text{supp }}
\def \ls{\lesssim}
\def\p{\partial}
\def\f{\frac}
\def\na{\nabla}
\def\alp{\alpha}

\def\O{\Omega}

\def\s{\sqrt}
\def\lab{\label}
\def\al{\aligned}
\def\eal{\endaligned}

 \allowdisplaybreaks

\begin{document}
 \footskip=0pt
 \footnotesep=2pt
\let\oldsection\section
\renewcommand\section{\setcounter{equation}{0}\oldsection}
\renewcommand\thesection{\arabic{section}}
\renewcommand\theequation{\thesection.\arabic{equation}}
\newtheorem{theorem}{\noindent Theorem}[section]
\newtheorem{lemma}{\noindent Lemma}[section]
\newtheorem*{claim}{\noindent Claim}
\newtheorem{proposition}{\noindent Proposition}[section]
\newtheorem*{inequality}{\noindent Weighted Inequality}
\newtheorem{definition}{\noindent Definition}[section]
\newtheorem{remark}{\noindent Remark}[section]
\newtheorem{corollary}{\noindent Corollary}[section]
\newtheorem{example}{\noindent Example}[section]

\title{Decay and vanishing of some  D-solutions of the Navier-Stokes equations}

\author{Bryan Carrillo \footnote{E-mail: bcarr011@ucr.edu},  Xinghong Pan\footnote{E-mail:xinghong{\_}87@nuaa.edu.cn},\q Qi S. Zhang\footnote{E-mail: qizhang@math.ucr.edu, corresponding address},\q Na Zhao\footnote{E-mail:nzhao13@fudan.edu.cn}\\
\vspace{0.5cm}
}

\date{}

\maketitle
% \vskip 0.2in

\centerline {\bf Abstract} An old problem since Leray \cite{Le:1} asks whether homogeneous
D-solutions of the 3 dimensional Navier-Stokes equation in $\bR^3$ or some noncompact domains are $0$. In this paper, we give a positive solution to the problem in two  cases: (1)  full 3 dimensional slab case $\bR^2 \times [0, 1]$ with Dirichlet boundary condition (Theorem \ref{ths2}); (2) when the solution is axially symmetric and periodic in the vertical variable (Theorem \ref{thzp}).

 Also, in the slab case, we prove that even if the Dirichlet integral has some growth, axially symmetric solutions with Dirichlet boundary condition must be swirl free, namely $u^\th=0$, thus reducing the problem to essentially a "2 dimensional" problem. In addition, a general D-solution (without the axial symmetry assumption) vanishes in $\R^3$ if, in spherical coordinates, the positive radial component of the velocity decays at order -1 of the distance.
 The paper is self contained comparing with  \cite{CPZ:1} although the general idea is related.

\vskip 0.3 true cm

\vskip 0.3 true cm

\vskip 0.3 true cm


\section{Introduction}
\q The purpose of the paper is to study decay and vanishing properties of the so called D-solutions to the steady Navier-Stokes equations
\be
\left\{
\begin{aligned}
&(u\cdot\nabla)u+\nabla p-\Delta u=f, \, \text{in} \, D \subset \bR^3 \\
&\nabla\cdot u=0,\\
\end{aligned}
\rt. \label{1.1}
\ee
with finite Dirichlet integral
\be
\int_{D}|\nabla u(x)|^2dx<+\infty \label{1.2}
\ee and various boundary conditions,  with also the requirement that $u$ vanishes at infinity; here $D \subset \bR^3$ is a noncompact or compact, connected  domain.  The two most basic noncompact domains are the whole space $\bR^3$ or a slab $\bR^2 \times [a, b]$; $f$ is a divergence free
 forcing term;
 Existence of these kind of solutions were studied in the pioneer work of Leray \cite{Le:1} (p24) by variational method and are often referred as D-solutions.
If $f=0$ and $u=0$ on $\p D$, then the solution is called a homogeneous D-solution. Given a noncompact domain, the following uniqueness problem has been open since then:

  \centerline{\it Is a homogeneous D-solution equal to  $0$ ?}

 This is also part of the very difficult uniqueness problem  for the steady Navier-Stokes equation.
In the 2 dimensional case, the problem in the full space case is solved by Gilbarg and Weinberger \cite{GW:1}. However, for the 3 dimensional problem, it is not even known if a general D-solution has any definite decay rate comparing with the distance function near infinity, even when the domain is $\bR^3$.
  The following is a list of  vanishing results with extra integral or decay assumptions for the solution $u$.  Galdi \cite{Ga:1} Theorem X.9.5 proved that if $u$ is a homogeneous D-solution in the domain $D=\bR^3$ and $u \in L^{9/2}(\bR^3)$, then $u=0$.  This result was improved by a log factor in Chae and Wolf \cite{CW:1}. In \cite{Ch:1}, Chae proved that homogeneous D-solutions in $\bR^3$ are $0$ if also $\Delta u \in L^{6/5}(\bR^3)$.  This condition scales the same way as $\Vert \nabla u \Vert_{L^2}$. Seregin \cite{Se:1} proved that homogeneous D-solutions in $\bR^3$ is $0$  if $u \in L^6(\bR^3) \cap BMO^{-1}$. In a recent paper \cite{KTW:1}, Kozono etc showed that homogeneous D-solutions in $\bR^3$ are $0$ if either the vorticity $w=w(x)$ decays faster than $c/|x|^{5/3}$ at infinity, or the velocity $u$ decays like $c/|x|^{2/3}$ with $c$ being a small number.
 Under certain smallness assumption, vanishing result for homogeneous 3 dimensional solutions in a slab was also obtained in the book \cite{Ga:1}, Chapter XII.

 The first result of the paper is a solution of the above problem if $D$ is a slab in $\R^3$.

\begin{theorem} \lab{ths2}
Let $u$ be a smooth, bounded solution to the problem
\be
\lab{3deqdiri}
\left\{
\begin{aligned}
&(u\cdot\nabla)u+\nabla p-\Delta u=0, \quad \text{in} \quad  \bR^2\times[0,1],\\
&\nabla\cdot u=0,\\
& u(x)|_{x_3=0}=u(x)|_{x_3=1}=0,
\end{aligned}
\rt.
\ee
such that the Dirichlet integral satisfies the condition:
\be\lab{d3d}
\int^1_{0}\int_{\bR^2}|\nabla u(x)|^2 dx<\infty.
\ee
Then, $u\equiv 0$.
\end{theorem}

{\remark Comparing with the full space case, one can show by Poincar\'e inequality that
the velocity $u$ is in $L^2$. Thus the decay rate of $u=u(x)$ is like $1/|x|$  in the integral sense. However one does not have a good knowledge of the pressure $p$. The Dirichlet boundary condition is known to induce complications on the vorticity and pressure. Our main work is to deal with the pressure term.  In addition, if the Dirichlet integral is infinite, then $0$ may not be the only solution. An example is $u=( x_3 (1-x_3), 0, 0)$, $p= -2 x_1$.
}

{\it The authors of the papers \cite{PS:1} and \cite{Pi:1} studied the asymptotic decay of solutions of the Navier-Stokes equation in a slab. They prove, under certain weighted integral assumption on the velocity $u=u(x)$ and 3rd order derivatives, $u=u(x)$ decays like $1/|x|$ or $1/|x|^3$ (See Theorem 3.1 in \cite{PS:1}). Then the vanishing of $u$ in the homogeneous case follows easily. However, these authors required that $(1+|x|)^{2+\beta} |u_3(x)|+(1+|x|)^{3+\beta} |\partial_{x_1} u_3(x)| $ with $\beta \in (-2, -1)$ is in $L^2$ in addition to further integral decay conditions of the first, second and third order derivatives of $u$, and consequently restriction on the pressure.  None of these conditions are available to us.
In the periodic case, which will be dealt later, it is not even known that $u$ is $L^2$.
}

Now if the domain is the whole $\bR^3$, we will show that if the positive part of the radial component of D-solutions decays at order $-1$ of the distance in spherical coordinates, then the D-solution vanishes.

\begin{theorem}\label{thgs}
Let $u_\rho=u_\rho(x)$ be the radial component of 3 dimensional D-solutions in spherical coordinates. If
\be\label{1.13}
{u_\rho(x)}\leq \f{C}{|x|}, \qquad x \in \bR^3,
\ee
for some positive constant $C$, then $u\equiv 0$.
\end{theorem}
\begin{remark}
 We should compare with the result in \cite{KTW:1} where the authors prove, if  the weak $L^{9/2}$  norm of $u$ is small, then $u$ vanishes. This includes the case $|u(x)| \le c |x|^{-2/3}$ for certain small constant $c$. In contrast,  our assumption here is worse on the order of the distance function. However we only impose the condition on the positive part of the radial component of the solution and there is no restriction on the other two components.
\end{remark}

 Next we concentrate on  axially symmetric homogeneous D-solutions, for which the vanishing problem is also wide open. Let us briefly describe the set up.
 It is convenient to work with cylindrical coordinates : $x=(x_1, x_2, z)$, $\th= \tan^{-1}{x_2/x_1}$  and $r=\sqrt{x^2_1+x^2_2}$. $e_r =(x_1/r, x_2/r, 0)$, $e_\th=(-x_2/r, x_1/r, 0)$ and $e_z=(0, 0, 1)$. We often write $x'=(x_1, x_2)$ and $x_3=z$.
A smooth vector field $u(x)=u^r(r,z)e_r+u^\th(r,z)e_\th+u^z(r,z)e_z$ is an axially symmetric solution of \eqref{1.1} if $u^r, u^\th, u^z$ satisfy the equations
\be \label{1.4}
\lt\{
\begin{aligned}
&(u^r\p_r+u^z\p_z)u^r-\f{(u^\th)^2}{r}+\p_r p=(\p^2_r+\f{1}{r}\p_r+\p^2_z-\f{1}{r^2})u^r,\\
&(u^r\p_r+u^z\p_z)u^\th+\f{u^ru^\th}{r}=(\p^2_r+\f{1}{r}\p_r+\p^2_z-\f{1}{r^2})u^\th,\\
&(u^r\p_r+u^z\p_z)u^z+\p_z p=(\p^2_r+\f{1}{r}\p_r+\p^2_z)u^z,\\
&\p_ru^r+\f{u^r}{r}+\p_zu^z=0.
\end{aligned}
\rt.
\ee Here and later we often refer to this equation as ASNS, i.e. axially symmetric Navier-Stokes equation.

The vorticity  $w$ is defined as $w(x)=\nabla\times u(x)=w^r(r,z)e_r+w^\th(r,z)e_\th+w^z(r,z)e_z$, where
\be\label{e1.4}
w^r=-\p_zu^\th,\q w^\th=\p_zu^r-\p_ru^z, \q w^z=\f{1}{r}\p_r(ru^\th).
\ee
The equations for $w^r, w^\th, w^z$ are
\be\label{1.5}
\lt\{
\begin{aligned}
&(u^r\p_r+u^z\p_z)w^r-(w^r\p_r+w^z\p_z)u^r=(\p^2_r+\f{1}{r}\p_r+\p^2_z-\f{1}{r^2})w^r,\\
&(u^r\p_r+u^z\p_z)w^\th-\f{u^r}{r}w^\th-\f{1}{r}\p_z(u^\th)^2=(\p^2_r+\f{1}{r}\p_r+\p^2_z-\f{1}{r^2})w^\th,\\
&(u^r\p_r+u^z\p_z)w^z-(w^r\p_r+w^z\p_z)u^z=(\p^2_r+\f{1}{r}\p_r+\p^2_z)w^z.
\end{aligned}
\rt.
\ee

 As mentioned earlier, although the components of solution $u$ is independent of the angle $\theta$ in the cylindrical system, the vanishing problem is also wide open. However certain decay estimates can be proven for the homogeneous D-solution $u$ and vorticity $w$.
 For example, the combined result of Chae-Jin \cite{CJ:1} and Weng \cite{We:1} state that, for $x \in \bR^3$,
\be
\lab{cjw}
|u(x)| \le C \left(\log r/r \right)^{1/2}, \q |w^\th(x)| \le C r^{-(19/16)^-}, \q
|w^r(x)| + |w^z(x)| \le C r^{-(67/64)^-}.
\ee Here $C$ is a positive constant and for a positive number $a$, we write $a^-$ as a number which is smaller than but close to $a$. Their proof is based on line integral techniques from Gilbarg and Weinberger \cite{GW:1}. In our previous work \cite{CPZ:1}, the decay estimate on $w$ is improved.  In that paper, we also proved a vanishing result when D-solutions are periodic in the third variable under the additional assumption that $u^\th$ and $u^z$ have zero mean in the $z$ direction. In a recent paper \cite{LRZ:1}, for solutions with infinite energy, Liouville property for  bounded,  axially symmetric solutions of the Navier-Stokes equation were studied under the natural assumption $r u^\theta$ is bounded.
Assuming in addition that $r u^\theta$ is bounded and $u$ is periodic in  $z$ variable, then it was shown that $u \equiv 0$. We mention that periodic solutions are also studied intensely in connection to the Kolmogorov flow. Earlier, the authors of the papers \cite{CSTY:1} and \cite{KNSS:1} proved Liouville theorems for ASNS under the assumption that $|u(x)| \le C/r$. See also an extension to $BMO^{-1}$ space in \cite{LZ:1}.

 The next result of this paper is an improvement of our main result in \cite{CPZ:1} in the case that the Dirichlet integral is finite. We will prove vanishing result when D-solutions are periodic in the third variable without any other assumption. Then we turn to axially symmetric solutions in a slab with Dirichlet boundary condition even if the Dirichlet integral has some growth. We will show that actually the  angular component of the solution vanishes.

 In the next theorem the flow is periodic in the $z$ direction with period $2 \pi$, a number chosen for convenience. Any other positive period also works. We will always take the forcing term $f=0$ throughout the paper.

\begin{theorem}\lab{thzp}
Let $u$ be a smooth axially symmetric solution to the problem
\be\label{1.7}
\left\{
\begin{aligned}
&(u\cdot\nabla)u+\nabla p-\Delta u=0, \quad \text{in} \quad  \bR^2 \times S^1=\bR^2 \times [-\pi, \pi],\\
&\nabla\cdot u=0,\\
&u(x_1,x_2,z)=u(x_1,x_2,z+2\pi),\\
&\lim\limits_{|x|\rightarrow \i}u=0,
\end{aligned}
\rt.
\ee
with finite Dirichlet integral:
\be
\int^\pi_{-\pi}\int_{\bR^2}|\nabla u(x)|^2dx<+\infty.\nn
\ee    Then $u=0$.
\end{theorem}
\begin{remark}
Assuming finiteness of the Dirichlet integral, Theorem\ref{thzp} is an improvement of our previous result Theorem 1.2 in \cite{CPZ:1} by removing the extra assumption that $u^\th$, $u^z$ have zero mean in the $z$ direction there. However,  the method in \cite{CPZ:1} is much different from the current one. Besides, the proof there can be applied to the case where the local Dirichlet integral has some growth, which means we can prove the vanishing result under the assumption $\int^\pi_{-\pi}\int_{|x'|\leq r}|\nabla u(x)|^2dx<(1+r)^\alp$ for some suitable and positive $\alp$. This method potentially allows application for flows with infinite Dirichlet energy such as Kolmogorov flows.
\end{remark}

The next theorem treats the case with Dirichlet boundary condition in a slab, even allowing the Dirichlet integral to  be log divergent.

\begin{theorem} \lab{thsd}
Let $u$ be a smooth, axially symmetric solution to the problem
\be\label{e1.8}
\left\{
\begin{aligned}
&(u\cdot\nabla)u+\nabla p-\Delta u=0, \quad \text{in} \quad  \bR^2\times[0,1],\\
&\nabla\cdot u=0,\\
&\lim\limits_{|x'|\rightarrow \i}u=0,\q u(x)|_{x_3=0}=u(x)|_{x_3=1}=0,
\end{aligned}
\rt.
\ee
such that the Dirichlet integral satisfies the condition: for a constant $C$, and all $R \ge 1$,
\be
\lab{dr2r}
\int^1_{0}\int_{R \le |x'| \le 2 R}|\nabla u(x)|^2 dx<C<\infty.
\ee
Then $u^\th = 0$. Moreover, there exists a positive constant $C_0$, depending only on the constant $C$ in (\ref{dr2r})  such that
\be
|u^r(x)| + |u^z(x)| \le C_0 \left(\f{\ln r}{r}\right)^{1/2}.
\ee
\end{theorem}
\begin{remark}
Since $u^\th = 0$ i.e. the flow is swirl free, the Navier-Stokes system reduces to
 \be\label{e1.9}
\lt\{
\begin{aligned}
&(u^r\p_r+u^z\p_z)u^r+\p_r p=(\Dl-\f{1}{r^2})u^r,\\
&(u^r\p_r+u^z\p_z)u^z+\p_z p=\Dl u^z,\\
&\p_ru^r+\f{u^r}{r}+\p_zu^z=0,\\
&\lim\limits_{r\rightarrow \i}(u^r,u^z)=0,\q (u^r,u^z)(r,z)|_{z=0,1}=0.
\end{aligned}
\rt.
\ee

So our vanishing problem now is much like a two dimensional problem. But unfortunately, we do not know any vanishing result for swirl free case in a slab with Dirichlet boundary condition and
(\ref{dr2r}).
\end{remark}

\begin{remark}
Clearly, if the Dirichlet integral is finite i.e. $\Vert \nabla u \Vert_{L^2(\bR^2 \times [0, 1])} < \infty$, then condition (\ref{dr2r}) is satisfied. If one works a little harder, then one can reach the same conclusion as the theorem assuming the integral in (\ref{dr2r}) grows at certain power of $R$. As mentioned earlier, solutions with infinite Dirichlet energy is also of interest. The decay estimate still holds if there is an inhomogeneous term of sufficiently fast decay.
\end{remark}

Now we outline the proof of the above results briefly. We start with the observation that in $z-$periodic case such that $D=\bR^2\times S^1$, the horizontal radial component of the solution $u^r$ satisfies $\int^\pi_{-\pi} u^rdz=0$. Poincar\'{e} inequality and the finite Dirichlet integral condition indicate that $u^r\in L^2(D)$. Then due to the speciality of $z$-periodic solution, we can actually prove that the oscillation of the pressure $p$ is bounded in a dyadic annulus. At last by testing the vector equation \eqref{1.7} with $u\phi^2(|x'|)$, where $\phi(x')$ is supported in $\{x'||x'|<2R\}$ and equal to $1$ in $\{x'||x'|<R\}$, and making $R$ approach  $\i$, we can prove that $u\equiv 0$.

In the case that $D=\bR^2\times [0,1]$ and $u$ with the homogeneous Dirichlet boundary condition,  we will show that the decay rate of $u^\th$ is $r^{-(\f{3}{2})^-}$ for large $r$. At the same time, the quantity $\G:=ru^\th$ satisfies
 \be\label{1.9}
(u^r\p_r+u^z\p_z)\G-(\Dl-\f{2}{r})\G=0.
\ee
$\G$ enjoys maximum principle, which means
\be\label{1.10}
\sup\limits_{x\in\O}|\G|\leq \sup\limits_{x\in\p\O}|\G|.
\ee
By using the above maximum principle, we can have $u^\th\equiv 0$.

Several steps are needed to get the decay of $u^\th$. In step one: the Green's function $G$ on $\bR^2\times [0,1]$ with homogeneous Dirichlet boundary condition will be introduced and a series of properties of $G$ will be displayed. The key point is that $G$ and its horizontal gradient have one more order decay on $|x'-y'|$ when $|x'-y'|>1$. In step two:  we obtain decay of $w^r,w^z$ by using a refined $Brezis-Gallouet\ inequality$, energy methods and scaling techniques to show that
\be
|(w^r,w^z)|\ls r^{-1}\ln r,\nn
\ee
 for large $r$. Furthermore, by the same procedure, we can show that
 \be
 |(\p_z w^r,\p_r w^z)|\ls r^{-3/2}(\ln r)^{3/2}.\\
 \ee
In step three, we use $Biot-Savart\ law$ to get the representation of $u^\th$ by $G$ and $(\p_z w^r,\p_r w^z)$ which implies that $u^\th$  decay in the same rate as $\p_z w^r$ and $\p_r w^z$.
Then we can apply the maximum principle on the function $\Gamma = r u^\th$ to conclude $u^\th=0$.

We conclude the introduction with a list of frequently used notations. For $(x, t) \in \bR^3 \times \bR$ and $r>0$, we use $Q_r(x, t)$ to denote the standard parabolic cube $\{ (y, s) \q | \q |x-y|<r, 0<t-s<r^2 \}$; $B(x, r)$ is the ball of radius $r$ centered at $x$. The symbol $ ... \ls ...$ stands for $... \le C ...$ for a positive constant $C$. $C$ with or without an index denotes a positive constant whose value may change from line to line. For $x=(x_1, x_2, x_3) \in \bR^3$, we write $x=(x', x_3)$ or $x=(x', z)$, and for $y=(y_1, y_2, y_3) \in \bR^3$, we write $y=(y', y_3)$. Theorem 1.1, 1.2, 1.3 and 1.4 will be proven in Sections 2, 3, 4 and 5 respectively.

\section{Proof of Theorem \ref{ths2}}

First, we will show two basic estimates for the $L^2$ norm of the velocity $u$ and the pressure $p$.
\subsection{ $\boldsymbol{L^2}$ estimates of $\boldsymbol{u}$ and $\boldsymbol{p}$}

Define the domain $\O_R=\big\{x'\in\bR^2\big||x'|\leq R\big\}\times[0,1]$, where $R\geq 1$. We consider the following problem:
 \\
\emph{Given}:
\be\label{e4.1e}
f\in L^2(\O_R) \qq {\rm with}\qq \int_{\O_R} f=0,
\ee
\emph{find a vector field $V:\O_R\rightarrow \bR^3$ such that}
\be\label{e4.2e}
\na\cdot V=f,\q V\in W^{1,2}_0(\O_R),\qq \|\na V\|_{L^2}\leq c_0 \|f\|_{L^2},
\ee
\emph{with} $c_0=c_0(\O_R)$.   For our purpose, we need an explicit estimate of the $c_0$ constant depending on the horizontal radius $R$.

The first solution of this problem is given in Bogovski\u{i} \cite{Bme:1,Bme:2}. See also Lemma III.3.1 of \cite{Ga:1}.
\begin{lemma}\label{L4.1L}
Let $\O=\big\{\bar{x}'\in\bR^2\big||\bar{x}'|\leq 1\big\}\times[0,1]$. Then for any $\bar{f}\in L^2(\O)$, satisfying
\be
\bar{f}\in L^2(\O) \qq {\rm with}\qq \int_{\O} \bar{f}=0,\nn
\ee
  there exists a constant $C$ and a vector valued function $\bar{V}:\O\rightarrow \bR^3$ such that
  \be\label{e4.3e}
\na\cdot \bar{V}=\bar{f},\q \bar{V}\in W^{1,2}_0(\O),\qq \|\na \bar{V}\|_{L^2}\leq C \|\bar{f}\|_{L^2},
\ee
where $C$ is an absolute constant.
\end{lemma}

Here we only stated a special case of Lemma III.3.1 of \cite{Ga:1} which works for general domains and the constant  $C$ depends on the ratio of the diameter and inner radius of the domain. For the above $\O$ the diameter and inner radius are fixed, so we have an absolute constant.  Note that the Bogovskii function may not be unique. For our purpose, we will choose and fix one for the relevant domain.

Next, we use the above lemma and a scaling argument to deduce the following proposition. We mention that the constant on the righthand side is $C R$ which  improves the one in Lemma III.3.1 of \cite{Ga:1}, which is $C R^4$. This is crucial for the proof of the vanishing result.

\begin{proposition}\label{P4.1}
Let $\O_R$ be as above. Then for any $f\in L^2(\O_R)$, satisfying \eqref{e4.1e}. Problem \eqref{e4.2e} has at least one solution $V$. Moreover, for the constant $c_0$ in \eqref{e4.2e}, we have the following estimate
\be\label{e4.4e}
\|\na V\|_{L^2(\O_R)}\leq  C R\|f\|_{L^2(\O_R)},
\ee
where $C$ is independent of $\O_R$.
\end{proposition}
\pf The existence of $V$ is already known as explained above. So we just need to prove (\ref{e4.4e}).

For $\bar{x}=(\bar{x}_1,\bar{x}_2,x_3)\in \O$, define
\[
\bar{f}(\bar{x}_1,\bar{x}_2,x_3):=f(R\bar{x}_1,R\bar{x}_2,x_3)= f(x_1, x_2, x_3).
\]Note $x_1=R\bar{x}_1$, $x_2=R\bar{x}_2$ but $x_3$ does not change.

It is easy to see that $\bar{f}$ satisfies the assumption in Lemma \ref{L4.1L}. So by Lemma \ref{L4.1L}, there exists a vector function $\bar{V}:\O\rightarrow \bR^3$ satisfying \eqref{e4.3e}. Then for $x\in\O_R$, define
\bea
&&V(x_1,x_2,x_3)=(V^1(x_1,x_2,x_3),V^2(x_1,x_2,x_3),V^3(x_1,x_2,x_3))\nn\\
&&\qq\qq\q\ \ =(R\bar{V}^1(\f{x_1}{R},\f{x_2}{R},x_3),R\bar{V}^2(\f{x_1}{R},\f{x_2}{R},x_3),\bar{V}^3(\f{x_1}{R},\f{x_1}{R},x_3)).
\eea
By a direct computation, we have
\be
\na\cdot {V}={f},\q {V}\in W^{1,2}_0(\O_R),\qquad \text{in} \quad {x} \quad \text{variables}\nn
\ee
\be
\na\cdot {\bar V}= \bar{f},\q \bar{V}\in W^{1,2}_0(\O_1), \qquad \text{in} \quad \bar{x} \quad \text{variables},\nn
\ee where $\bar{V} = (\bar{V}^1(\bar{x}), \bar{V}^2(\bar{x}), \bar{V}^3(\bar{x}))$.
Now we  estimate the $L^2$ norm of $\na V$. We use $\alpha,\beta$ to take values only on $1,2$ and $i,j$ to take values on $1,2,3$. So we have
\bea\label{e4.6e}
&&\|\na V\|^2_{L^2(\O_R)}=\sum\limits^3_{i,j=1}\int^1_0\int_{|x'|\leq R}|\frac{\p V^j}{\p x_i}|^2dx'dx_3\nn\\
&&\qq\qq\q\ \ =\int^1_0\int_{|x'|\leq R}\big(\sum\limits^2_{\alpha,\beta=1}| \frac{\p V^\beta}{\p x_\alpha}|^2+\sum\limits^2_{\beta=1}|\frac{\p V^\beta}{\p x_3}|^2\nn\\
&&\qq\qq\qq\qq\qq+\sum\limits^2_{\alpha=1}|\frac{\p V^3}{\p x_\alpha}|^2dx'dx_3+|\frac{\p V^3}{\p x_3}|^2\big)dx'dx_3\nn\\
&&\qq\qq\q\ \ =\int^1_0\int_{|x'|\leq R}\Big(\sum\limits^2_{\alpha,\beta=1}|\frac{\p \bar{V}^\beta}{\p \bar{x}_\alpha}|^2({\f{x'}{R}},x_3)+R^2\sum\limits^2_{\beta=1}|\frac{\p \bar{V}^\beta}{\p x^3}|^2({\f{x'}{R}},x_3)\nn\\
&&\qq\qq\qq\q+\f{1}{R^2}\sum\limits^2_{\alpha=1}|\frac{\p \bar{V}^3}{\p \bar{x}^\alpha}|^2({\f{x'}{R}},x_3)+|\frac{\p \bar{V}^3}{\p x_3}|^2({\f{x'}{R}},x_3)\Big)dx'dx_3\nn\\
&&\qq\qq\q\ \ =R^2\int^1_0\int_{|\bar{x}'|\leq 1}\Big(\sum\limits^2_{\alpha,\beta=1}|\frac{\p \bar{V}^\beta}{\p \bar{x}^\alpha}|^2(\bar{x}',x_3)+R^2\sum\limits^2_{\beta=1}|\frac{\p \bar{V}^\beta}{\p x_3}|^2(\bar{x}',x_3)\nn\\
&&\qq\qq\qq\q+\f{1}{R^2}\sum\limits^2_{\alpha=1}|\frac{\p \bar{V}^3}{\p \bar{x}^\alpha}|^2(\bar{x}',x_3)+|\frac{\p \bar{V}^3}{\p x_3}|^2(\bar{x}',x_3)\Big)d\bar{x}'dx_3\nn\\
&&\qq\qq\qq\leq CR^4\|\na \bar{V}\|^2_{L^2(\O_1)}.
\eea
Also it is easy to see
\be\label{e4.7e}
\|f\|^2_{L^2(\O_R)}=R^2\|\bar{f}\|^2_{L^2(\O_1)}.
\ee
Combining \eqref{e4.6e}, \eqref{e4.7e} and \eqref{e4.3e}, we have
\be
\|\na V\|^2_{L^2(\O_R)}\leq CR^4\|\na \bar{V}\|^2_{L^2(\O)}\leq CR^4\|\bar{f}\|^2_{L^2(\O)}=CR^{2}\|f\|^2_{L^2(\O_R)}.\nn
\ee
This finishes the proof of Proposition \ref{P4.1}. \qed

Next we give some $L^2$ estimate of the velocity $u$ and the pressure $p$, using the preceding proposition.

\begin{lemma}
\lab{leupl3}
Let $u,p$ be the solution of \eqref{3deqdiri}, then we have
\be
\lab{ul2}
\int_{\bR^2\times[0,1]} |u|^2dx<\infty,
\ee

\be\label{e4.9e}
\|p-p_R\|_{L^2(\O_R)}\leq C_0R,
\ee
where $C_0=C(\|u\|_{L^\i},\|u\|_{L^2},\|\na u\|_{L^2})$ and $p_R:=\f{1}{|\O_R|}\int_{\O_R}pdx$ is the average of $p$ on $\O_R$.
\end{lemma}
\pf
Since we have zero boundary on $x_3=0,1$, the one dimensional Poincar\'e inequality indicates that
\bea
&&\int_{\bR^2\times[0,1]} |u|^2dx=\int_{\bR^2}\int^1_0 |u|^2dx_3dx'\nn\\
&&\qq\qq\qq\ls\int_{\bR^2}\int^1_0 |\p_z u|^2dx_3dx'\nn\\
&&\qq\qq\qq\ls \int_{\bR^2\times[0,1]} |\na u|^2dx<\infty
\eea by the definition of D-solutions. This proves (\ref{ul2}).

From Proposition \ref{P4.1}, there exists a $V$ satisfying \eqref{e4.2e} and \eqref{e4.4e} with $f=p-p_R$. Now multiplying $\eqref{e1.8}_1$ with $V$ and integration on $\O_R$, we get
\be
\int_{\O_R}\na(p-p_R)\cdot V dx=\int_{\O_R}(\Dl u-u\cdot\na u)\cdot V dx.
\ee
Integration by parts indicate that
\bea
\int_{\O_R}(p-p_R)^2dx&=&\int_{\O_R}(p-p_R)\na\cdot V dx\nn\\
&=&\int_{\O_R}\sum^3_{i,j=1}\p_iu^j\p_iV^j+\na\cdot(u\otimes u)\cdot V dx\nn\\
&=& \int_{\O_R}\sum^3_{i,j=1}(\p_iu^j-u^iu^j)\p_iV^j dx\nn\\
&\leq&\|\na V\|_{L^2(\O_R)}\big(\|\na u\|_{L^2(\O_R)}+\|u\|_{L^\i(\O_R)}\| u\|_{L^2(\O_R)}\big)\nn\\
&\leq& \f{\ve}{R^2}\|\na V\|^2_{L^2(\O_R)}+C_\ve R^2\big(\|\na u\|_{L^2(\O_R)}+\|u\|_{L^\i(\O_R)}\| u\|_{L^2(\O_R)}\big)^2\nn\\
&\leq& C\ve\|p-p_R\|^2_{L^2(\O_R)}+C_\ve R^2\big(\|\na u\|_{L^2(\O_R)}+\|u\|_{L^\i(\O_R)}\| u\|_{L^2(\O_R)}\big)^2.\nn
\eea Here to reach the last line, we used \eqref{e4.4e}. By choosing $\ve$ is small enough, we can obtain \eqref{e4.9e}.\ef

\subsection{Vanishing of $\boldsymbol{u}$}

Now we are in a position to complete the proof of Theorem \ref{ths2}. Let $\phi(s)$ be a smooth cut-off function satisfying
\be\label{testf1}
\phi(s)=\lt\{
\begin{aligned}
& 1\qq s\in [0,1/2],\\
& 0 \qq   s\geq 1,
\end{aligned}
\rt.
\ee with the usual property that $\phi$, $ \phi'$ and $\phi''$ are bounded.
Set $\phi_R(y')=\phi(\f{|y'|}{R})$ where $R$ is a large positive number. For convenience of notation, we denote $I=[0,1]$. Now testing the Navier-Stokes equation
\[
u\cdot\na u+\na p=\Dl u
\]with $u\phi_R$, we obtain
\be
\int_{\bR^2\times I}-\Dl u (u\phi_R)dx=\int_{\bR^2\times I}-(u\cdot\na u+\na (p-p_R) (u\phi_R)dx. \nn
\ee
Integration by parts indicates that
\be
\begin{aligned}
&\q\int_{\bR^2\times I}|\na u|^2 \phi_R dx-\f{1}{2}\int_{\bR^2\times I}|u|^2\Dl\phi_R dx\\
&=-\f{1}{2}\int_{\bR^2\times I}u\cdot\na |u|^2\phi_R dx+\int_{\bR^2\times I} (p-p_R)u\cdot\na\phi_Rdx\\
&=\f{1}{2}\int_{\bR^2\times I} |u|^2u\cdot\na\phi_R dx+\int_{\bR^2\times I} (p-p_R)u\cdot\na \phi_R dx.
\end{aligned}\nn
\ee
Denote $\bar{B}_{R/\f{1}{2}R}:=\{x'|1/2R\leq|x'|\leq R\}$ i.e. the dyadic annulus. Then we have, since
$\phi_R$ depends only on $r$, that
\be
\begin{aligned}
&\q\int_{\bR^2\times I}|\na u|^2 \phi_Rdx\\
&\ls \f{1}{R^2}\int^{1}_{0}\int_{\bar{B}_{R/\f{1}{2}R}}|u|^2dx+
\f{1}{R}\int^{1}_{0}\int_{\bar{B}_{R/\f{1}{2}R}} |u^r| \, |u|^2dx\\
&+ \f{1}{R}\int^{1}_{0}\int_{\bar{B}_{R/\f{1}{2}R}}|p-p_R||u^r|dx\\
&\ls \f{1}{R^2} \int^{1}_{0}\int_{\bar{B}_{R/\f{1}{2}R}} |u|^2 dx +
\f{\Vert u^r \Vert_\infty}{R}
\int^{1}_{0}\int_{\bar{B}_{R/\f{1}{2}R}} |u|^2 dx  \\
&\q+
\f{1}{R} \left(\int^{1}_{0}\int_{\bar{B}_{R/\f{1}{2}R}} (u^r)^2 dx \right)^{1/2}
\left(\int^{1}_{0}\int_{\bar{B}_{R/\f{1}{2}R}}|p-p_R|^2dx \right)^{1/2}.\\
\end{aligned}\nn
\ee
By the estimate of pressure (\ref{e4.9e}), we deduce:
\be
\begin{aligned}
\q\int_{\bR^2\times I}|\na u|^2 \phi_R dx &\ls \big(\Vert u \Vert^2_{L^\infty(\bar{B}_{R/\f{1}{2}R}\times I)}+C_0\big)\Vert u \Vert_{L^2(\bar{B}_{R/\f{1}{2}R}\times I)}.
\end{aligned}\nn
\ee
Now let $R\rightarrow +\i$, using $u\in L^2(\bR^2\times I)$ (Lemma \ref{leupl3}), we arrive at
\be
\int_{\bR^2\times I}|\na u|^2dx=0,\nn
\ee
which shows that $u\equiv c.$ Besides, recall $u=0$ at the boundary, then at last we deduce
\be
u\equiv 0.
\ee This completes the proof of Theorem \ref{ths2}.
\ef

\section{Proof of Theorem \ref{thgs}}

\q In this section, we prove Theorem\ref{thgs}. Recall the following result proved by Galdi (see Theorem X.5.1 of \cite{Ga:1} for a more general version).
\begin{proposition}
Let $u(x)$ be a generalized solution of (NS) satisfying \eqref{1.1} and $p(x)$ be the associated pressure, then there exists $p_1\in\bR$ such that
\be
\lim\limits_{|x|\rightarrow \i} |D^\alp u(x)| + \lim\limits_{|x|\rightarrow \i}|D^\alp (p(x)-p_1)|=0
\ee
uniformly for all multi-index $\alp = (\alp_1, \alp_2, \alp_3)\in [\bN\cup\{0\}]^3$.
\end{proposition}

Here the notion of generalized solution is the one given in Definition X.1.1 \cite{Ga:1} p653.

Define the head pressure $Q:=\f{1}{2}|u|^2+p-p_1$, which satisfies
\be\label{3.3}
-\Dl Q+ u \cdot\na Q=-|\text{curl} \, u|^2.
\ee
From the preceding Proposition, we have $\lim\limits_{|x|\rightarrow \i} Q=0.$ \eqref{3.3} indicates that $Q$ satisfies the maximum principle. So we have that
\be\label{3.4}
Q\leq 0.
\ee
Set $\phi_R(x)=\phi(\f{|x|}{R})$, where $\phi$ is defined in \eqref{testf1}. Now testing the Navier-Stokes equation $u\cdot\na u+\na p=\Dl u$ with $u\phi_R$, we can get
\be
\int_{\bR^3}-\Dl u (u\phi_R)dx=\int_{\bR^3}-(u\cdot\na u+\na p) (u\phi_R)dx.\nn
\ee
Integration by parts indicates that
\be
\begin{aligned}
&\q\int_{\bR^3}|\na u|^2 \phi_Rdx-\f{1}{2}\int_{\bR^3}|u|^2\Dl\phi_R dx\\
&=-\f{1}{2}\int_{\bR^3}u\cdot\na |u|^2\phi_Rdx+\int_{\bR^3} (p-p_1)u\cdot\na\phi_Rdx\\
&=\f{1}{2}\int_{\bR^3} |u|^2u\cdot\na\phi_Rdx+\int_{\bR^3} (p-p_1)u\cdot\na\phi_Rdx\\
&=\int_{\bR^3} Q u\cdot\na\phi_Rdx=\int_{\bR^3} Q u_\rho\p_\rho \phi_R dx.
\end{aligned}\nn
\ee
Then we get
\be\label{3.5}
\begin{aligned}
&\q\int_{\bR^3}|\na u|^2 \phi_Rdx\\
&\ls \f{1}{R^2}\int_{B_{R/\f{1}{2}R}}|u|^2 dx - \int_{\bR^3} Q u^-_\rho\p_\rho \phi_R dx+\int_{\bR^3} Q u^+_\rho\p_\rho \phi_R dx,
\end{aligned}
\ee
where $u^-_\rho=:-\min\{0,u_\rho\}$ and $u^+_\rho=:\max\{0,u_\rho\}$. Also
$B_{R/\f{1}{2}R} \equiv B_{R}-B_{\f{1}{2}R}$ here and later in the proof.

Since $u\in L^6(\bR^3),\ p-p_1\in L^3(\bR^3)$, we have $Q\in L^3(\bR^3)$. Besides, it is obvious that we can choose $\p_\rho \phi_{R}\leq 0$. So by using \eqref{3.4} and from \eqref{3.5}, we get
\be\label{3.6}
\begin{aligned}
&\q\int_{\bR^3}|\na u|^2 \phi_Rdx\\
&\ls \f{1}{R^2}\Big(\int_{B_{R/\f{1}{2}R}}|u|^6 dx\Big)^{1/3}\Big(\int_{B_{R/\f{1}{2}R}}dx\Big)^{2/3} - \int_{\bR^3} Q u^-_\rho\p_\rho \phi_R dx\\
&\q +\f{1}{R}\sup\limits_{B_{R/\f{1}{2}R}}{u^+_\rho}\Big(\int_{B_{R/\f{1}{2}R}} |Q|^3 dx\Big)^{1/3}\Big(\int_{B_{R/\f{1}{2}R}} dx\Big)^{2/3}\\
& \ls \Big(\int_{B_{R/\f{1}{2}R}}|u|^6 dx\Big)^{1/3}+\Big(\int_{B_{R/\f{1}{2}R}} |Q|^3 dx\Big)^{1/3}.
\end{aligned}
\ee
Now let $R\rightarrow +\i$, we have
\be
\int_{\bR^3}|\na u|^2dx=0,\nn
\ee
which implies that $u\equiv c.$ Since $\lim\limits_{|x'|\rightarrow\i}u=0$, then at last we deduce
\be
u\equiv 0.
\ee
\ef

\section{Proof of Theorem \ref{thzp}} \label{s2}

 This section is divided into 2 parts.  First, by integrating  the first and third equation of \eqref{1.4}, we can prove that actually the oscillation of the pressure $p$ is bounded in dyadic annulus in the axially symmetric $z-$periodic case. Second by testing the Navier Stokes equations by a standard test function,  we can prove the vanishing result.

 Let us remark that boundedness of $u$, $\na u$ and $\na^2 u$ are known, as given  in the interior regularity theorem of \cite{Ga:1} (Theorem  XI.1.2, p755).
  Alternatively one can apply a result for local solutions in \cite{Z:2}. The point is that no additional information is needed for the pressure $p$.

\subsection{Boundedness of the oscillation of $\boldsymbol{p}$ in dyadic annulus}

In this subsection we prove that the oscillation of the pressure $p$ is bounded in dyadic annulus $B(0, 2R)-B(0, R)$. In fact, if one uses the a priori decay of solutions in Corollary 2.1 of \cite{CPZ:1}, one can show that the oscillation of $p$ in the whole domain is bounded. However in order to keep the proof independent and short, we will not prove that here.

From the third equation of \eqref{1.4} and using the boundedness of $u,\na u$ and $\na^2 u$ mentioned above, we have
\be\label{2.20'}
|\p_z p|\ls 1.
\ee
Next we will show that for any $R>1$, $R<r<2R$, we have
\be\label{e2.12}
\Big|\int^{\pi}_{-\pi}(p(r,z)-p(R,z))dz\Big|\ls 1.
\ee
Integrating the first equation of \eqref{1.4} on $z$ from $-\pi$ to $\pi$, we can get
\be\label{ee2.21}
\begin{aligned}
&\p_r\int^\pi_{-\pi}pdz=\int^\pi_{-\pi} \left[ -(u^r\p_r+u^z\p_z)u^r+\f{(u^\th)^2}{r}+(\p^2_r+\f{1}{r}\p_r+\p^2_z-\f{1}{r^2})u^r \right] dz\\
&\qq\qq\ =-\int^\pi_{-\pi}\f{1}{2}\p_r (u^r)^2dz-\int^\pi_{-\pi}u^z\p_z u^rdz+\int^\pi_{-\pi}\f{(u^\th)^2}{r}dz\\
&\qq\qq\q+\p^2_r\int^\pi_{-\pi}u^rdz+
\f{1}{r}\p_r\int^\pi_{-\pi}u^rdz-\f{1}{r^2}\int^\pi_{-\pi}u^rdz.
\end{aligned}
\ee
Picking any $r_0 \in [R, 2 R]$, integrating \eqref{ee2.21} on $r$ from $R$ to $r_0$, we find
\be\label{ee2.22}
\begin{aligned}
&\q\big|\int^\pi_{-\pi}p(r_0,z)dz-\int^\pi_{-\pi}p(R,z)dz\big|\\
& \ls\underbrace{\big|\int^{r_0}_R\int^\pi_{-\pi}\p_r (u^r)^2dzdr\big|}_{I_1}+\underbrace{\big|\int^{r_0}_R\int^\pi_{-\pi}u^z\p_z u^rdzdr\big|}_{I_2}\\
&\q+\underbrace{\big|\int^{r_0}_R\int^\pi_{-\pi}\f{(u^\th)^2}{r}dzdr\big|}_{I_3}
+\underbrace{\big|\int^{r_0}_R\p^2_r\int^\pi_{-\pi}u^rdzdr\big|}_{I_4}\\
&\q+\underbrace{\big|\int^{r_0}_R\f{1}{r}\p_r\int^\pi_{-\pi}u^r dzdr\big|}_{I_5}+\underbrace{\big|\int^{r_0}_R\f{1}{r^2}\int^\pi_{-\pi}u^rdzdr\big|}_{I_6}.
\end{aligned}
\ee
Next we show that all the terms on the right hand of \eqref{ee2.22} are bounded.\\
Changing the integration order of $r$ and $z$, we have
\be
\begin{aligned}
I_1=\big|\int^\pi_{-\pi}\big[(u^r)^2(r_0, z)- (u^r)^2(R,z)\big]dz \big|\ls 1.
\end{aligned}
\ee
Using integration by parts and the incompressible condition, we can obtain
\be
\begin{aligned}
&I_2=\big|\int^{r_0}_R\int^\pi_{-\pi} u^r\p_z u^zdz dr\big|\\
&\q =\big|\int^{r_0}_R\int^\pi_{-\pi} u^r(\f{u^r}{r}+\p_r u^r)dz dr\big|\\
&\q \ls \big|\int^{r_0}_R\int^\pi_{-\pi} \f{(u^r)^2}{r}dz dr\big|+\big|\int^{r_0}_R\int^\pi_{-\pi} \p_r (u^r)^2dz dr\big|\\
&\q \ls  \big|\int^{2R}_R\int^\pi_{-\pi}r^{-1} drdz\big|+\big|\int^\pi_{-\pi}\big[(u^r)^2(r_0,z)- (u^r)^2(R,z)\big]dz \big|\\
&\q  \ls  1.
\end{aligned}
\ee
Using the boundedness of $u$, we deduce
\be
\begin{aligned}
&I_3 \ls \big|\int^{2R}_R\int^\pi_{-\pi}r^{-1} drdz\big| \ls  1.
\end{aligned}
\ee
Using the boundedness of $\na u$, we have
\be
\begin{aligned}
&I_4=\big|\int^\pi_{-\pi}[ \p_r u^r(r_0,z)- \p_ru^r(R,z)] dz \big|\ls 1.
\end{aligned}
\ee
Integration by parts implies that
\be
\begin{aligned}
&I_5=\big|\int^\pi_{-\pi} \big(\f{u^r}{r}\bigg|^{r_0}_R+\int^{r_0}_R\f{u^r}{r^2}dr\big)dz\big|\ls 1.
\end{aligned}
\ee
Also the boundedness of $u$ indicate that
\be
\begin{aligned}
&I_6\ls\big|\int^\pi_{-\pi} \int^{r_0}_1 r^{-2}drdz\big|\ls 1.
\end{aligned}
\ee

Combining the above, we obtain the boundedness of the right hand of \eqref{ee2.22},  which shows, for any $R>1, R<r<2R$, that
\be\label{ee2.29}
\Big|\int^\pi_{-\pi}(p(r,z)-p(R,z))dz\Big|\ls 1.
\ee
Then, according to the mean value theorem,
fixing $R>1$, for any $R<r<2R$, there exists $z(r)$, such that
\be\label{e2.30}
\big|p(r,z(r))-p(R,z(r))\big|\ls 1.
\ee

Combination of \eqref{e2.30} and uniform boundedness of $\p_z p$, we can get for $R>1$, any $R<r<2R$, and $z\in [-\pi,\pi]$, the following bound on oscillation of $p$:
\be\label{2.32}
\begin{aligned}
&|p(r,z)-p(R,0)|\\
&\qq\qq =|(p(r,z)-p(r,z(r)))+(p(r,z(r))-p(R,z(r)))\\
&\qq\qq\qq\qq + (p(R,z(r))-p(R, 0))|\\
&\qq\qq \leq (|\p_z p(r,z_1)|+|\p_z p(R,z_2)|)(|z-z(r)|+|z(r)|)+|p(r,z(r))-p(R,z(r))|\\
&\qq\qq\ls 1,
\end{aligned}
\ee
where $z_1,z_2\in [-\pi, \pi]$ and we have used the mean value theorem.
\subsection{Vanishing of $\boldsymbol{u}$}

First we use the one dimensional Poincar\'e inequality  and the boundedness of the Dirichlet integral to prove that
\be\label{2.33}
\int^\pi_{-\pi}\int_{\bR^2}|u^r|^2dx\ls 1.
\ee

To start, we make a claim that
\[
\int^\pi_{-\pi}u^r (r, z) dz =0, \quad \forall r \ge 0.
\]This can be seen from integrating the divergence free condition in the $z$ direction:
\[
\partial_r u^r + \frac{u^r}{r} + \partial_z u^z =0,
\]giving us, since $u^z$ is periodic in $z$:
\[
r \partial_r \int^\pi_{-\pi} u^r dz + \int^\pi_{-\pi} u^r dz =0.
\]Therefore
\[
\left( r \int^\pi_{-\pi} u^r(r, z) dz \right)'=0,
\]i.e.
\[
r \int^\pi_{-\pi} u^r(r, z) dz = 0 \times \int^\pi_{-\pi} u^r(0, z) dz =0.
\]This proves the claim.
 Alternately we can use $\int^\pi_{-\pi}u^r dz=-\int^\pi_{-\pi} \p_z L^\th dz=0$, where $L^\th$ is the angular stream function. But this requires to construct a global stream function through the equation
 \[
 \Delta L^\theta - \frac{1}{r^2} L^\theta = -w^\theta,
 \]which may take longer.

 Hence  we have
\be
\begin{aligned}
&\q\int^\pi_{-\pi}\int_{\bR^2}|u^r|^2dx\\
&=\int_{\bR^2}\int^\pi_{-\pi}\big(u^r-\int^\pi_{-\pi}u^rdz\big)^2dx'dz\\
&\ls \int_{\bR^2}\int^\pi_{-\pi}\big|\p_zu^r|^2dzdx'\\
&\ls \int^\pi_{-\pi}\int_{\bR^2}\big|\na u|^2dx\ls 1.
\end{aligned}\nn
\ee

Let $\phi(s)$ be a smooth cut-off function satisfying

\be\label{testf}
\phi(s)=\lt\{
\begin{aligned}
& 1\qq s\in [0,1],\\
& 0 \qq   s\geq 2,
\end{aligned}
\rt.
\ee with the usual property that $\phi$, $ \phi'$ and $\phi''$ are bounded.
Set $\phi_R(y')=\phi(\f{|y'|}{R})$ where $R$ is a large positive number. For convenience of notation, we denote $I=[-\pi,\pi]$. Now testing the Navier-Stokes equation
\[
u\cdot\na u+\na p=\Dl u
\]with $u\phi_R$, we obtain
\be
\int_{\bR^2\times I}-\Dl u (u\phi_R)dx=\int_{\bR^2\times I}-(u\cdot\na u+\na (p-p(R,0)) (u\phi_R)dx. \nn
\ee
Integration by parts indicates that
\be
\begin{aligned}
&\q\int_{\bR^2\times I}|\na u|^2 \phi_Rdx-\f{1}{2}\int_{\bR^2\times I}|u|^2\Dl\phi_R dx\\
&=-\f{1}{2}\int_{\bR^2\times I}u\cdot\na |u|^2\phi_R dx+\int_{\bR^2\times I} (p-p(R,0))u\cdot\na\phi_Rdx\\
&=\f{1}{2}\int_{\bR^2\times I} |u|^2u\cdot\na\phi_R dx+\int_{\bR^2\times I} (p(r,z)-p(R,0))u\cdot\na \phi_R dx.
\end{aligned}\nn
\ee
Denote $\bar{B}_{2R/R}:=\{x'|R\leq|x'|\leq 2R\}$ i.e. the dyadic annulus. Then we have, since
$\phi_R$ depends only on $r$, that
\be
\begin{aligned}
&\q\int_{\bR^2\times I}|\na u|^2 \phi_Rdx\\
&\ls \f{1}{R^2}\int^{\pi}_{-\pi}\int_{\bar{B}_{2R/R}}|u|^2dx+
\f{1}{R}\int^{\pi}_{-\pi}\int_{\bar{B}_{2R/R}} |u^r| \, |u|^2dx\\
&+\sup\limits_{R<r<2R,z\in[-\pi,\pi]}\big|(p(r,z)-p(R,0))
\big|\int^{\pi}_{-\pi}\int_{\bar{B}_{2R/R}}|u^r| |\p_r\phi_R| dx\\
&\ls \f{\Vert u \Vert^2_{L^\infty(\bar{B}_{2R/R} \times [-\pi, \pi])}}{R^2} \int^{\pi}_{-\pi}\int_{\bar{B}_{2R/R}}dx
 \\
&\q +
\f{\Vert u \Vert^2_{L^\infty(\bar{B}_{2R/R}\times [-\pi, \pi])}}{R} \left(\int^{\pi}_{-\pi}\int_{\bar{B}_{2R/R}} (u^r)^2 dx \right)^{1/2}
\left(\int^{\pi}_{-\pi}\int_{\bar{B}_{2R/R}}dx \right)^{1/2} \\
&\q+ \f{C_0}{R}\Big(\int^{\pi}_{-\pi}\int_{\bar{B}_{2R/R}}|u^r|^2dx\Big)^{1/2}
\Big(\int^{\pi}_{-\pi}\int_{\bar{B}_{2R/R}}dx\Big)^{1/2} \quad (\text{by} \quad (\ref{2.32}))\\
&\ls \Vert u \Vert^2_{L^\infty(\bar{B}_{2R/R} \times [-\pi, \pi])}
+C_0\Big(\int^{\pi}_{-\pi}\int_{\bar{B}_{2R/R}}|u^r|^2dx\Big)^{1/2}.
\end{aligned}\nn
\ee Here
\[
C_0 = C \sup\limits_{R<r<2R,z\in[-\pi,\pi]}\big|(p(r,z)-p(R,0))| + C \Vert u \Vert^2_{L^\infty(\bar{B}_{2R/R}\times[-\pi,\pi])}.
\]
Now let $R\rightarrow +\i$, using (\ref{2.33}) and the assumption that $u \to 0$ as $r \to \infty$, we arrive at
\be
\int_{\bR^2\times I}|\na u|^2dx=0,\nn
\ee
which shows that $u\equiv c.$ Besides, we have $\lim\limits_{|x'|\rightarrow\i}u=0$, then at last we deduce
\be
u\equiv 0.
\ee This completes the proof of the theorem.
\ef

\section{Proof of Theorem \ref{thsd}}
\subsection{The Green function in $\bR^2\times [0,1]$}

\begin{lemma}\label{LemmaG}
 Let $G=G(x,y)$ be the distribution solution of the following equation, namely the Green's function.
\be\label{GFD}
\lt\{
\begin{aligned}
&-\Dl G(x,y)=\dl(x,y),\qq y\in\bR^2\times[0,1],\\
&G(x,y)|_{y_3=0}=G(x,y)|_{y_3=1}=0
\end{aligned}
\rt.
\ee

Then $G(x,y)$ has the following representation formula
\be
G(x,y)
=\f{1}{4\pi}\sum\limits^{+\i}_{n=-\i}\bigg\{\f{1}{\s{|x'-y'|^2+|x_3-y_3+2n|^2}}
-\f{1}{\s{|x'-y'|^2+|x_3+y_3-2n|^2}}\bigg\}.\nn
\ee

\end{lemma}

\pf We can use Method of Images to deduce the precise formula of $G(x,y)$. Readers can see Example 2.13 in \cite{MM:1} for a proof in $\bR^1\times[0,1]$. The derivation for the case $\bR^2\times [0,1]$ is similar, we omit the details. \ef

\begin{lemma} \label{PG}
Let $G(x,y)$ be the Green function in Lemma\ref{LemmaG}, then there exists a $c>0$ such that
\be\label{4.3'}
|G(x,y)| \le \f{c}{|x'-y'|(1+|x'-y'|)+|x_3-y_3|}
\ee
and
\be\label{4.4}
|\p_{x',y'}G(x,y)| \le \f{c}{|x'-y'|^2(1+|x'-y'|)+|x_3-y_3|^2}
\ee
with $x'=(x_1,x_2)$ and $y'=(y_1,y_2)$.
\end{lemma}

\pf Since the proof is a series of tedious computation, we put it in the \textbf{Appendix}.
Faster decay for the Green's function is also true, but we will not need it here. \qed

\begin{lemma}\label{l10.3}
Let ${G}(x,y)$ be defined as above. Denote $x=(r\cos\th,r\sin\th,z)$ and $y=(\rho\cos\phi,\rho\sin\phi,l)$, then we have the following estimates
when $|\rho-r|\leq \f{1}{4}r$:

\be\label{4.6}
\int^{2\pi}_0|{G}(x,y)|d\phi\ls \lt\{
\begin{aligned}
&\f{1}{r|\rho-r|} \qq\qq\q\  1\leq|\rho-r|\leq \f{r}{4},\\
&\f{1}{r}\ln (1+\f{r}{|\rho-r|})\qq  |\rho-r|<1,
\end{aligned}
\rt.
\ee
\be\label{4.7}
\int^{2\pi}_0|\p_{\rho,r} {G}(x,y)|d\phi\ls \lt\{
\begin{aligned}
&\f{1}{r|\rho-r|^2} \qq\qq\qq 1\leq|\rho-r|\leq \f{r}{4},\\
&\f{1}{r(|\rho-r|+|z-l|)}\qq  |\rho-r|<1.
\end{aligned}
\rt.
\ee
\end{lemma}
\pf
Remember that $x'=(r\cos\th,r\sin\th), y'=(\rho\cos\phi, \rho \sin\phi)$, then we can get
\be
|x'-y'|=\s{\rho^2+r^2-2\rho r\cos(\th-\phi)}\xlongequal{\th=0}\s{(\rho-r)^2+4\rho r\sin^2\f{\phi}{2}}.\nn
\ee
When $|\rho-r|<1$, from Lemma \ref{PG}, it is easy to see that
\be\label{10.19}
\begin{aligned}
&\q \int^{2\pi}_0|{G}(x,y)|d\phi
\ls \int^{2\pi}_0\f{1}{|x'-y'|}d\phi\\
&\ls\int^{\pi}_0 \f{d\phi}{\s{(\rho-r)^2+4\rho r\sin^2\phi}}\\
&\ls\int^{\f{\pi}{2}}_0 \f{d\phi}{\s{(\rho-r)^2+4\rho r\sin^2\phi}}\\
& \thickapprox \f{1}{r}\int^{\f{\pi}{2}}_0 \f{d\phi}{\s{\kappa^2+\sin^2\phi}} \q \text{with}\q \kappa^2=\f{|\rho-r|^2}{4\rho r}\ls 1\\
& \thickapprox \f{1}{r}\Big(\int^{\f{\pi}{4}}_0+\int^{\f{\pi}{2}}_{\f{\pi}{4}}\Big) \f{d\phi}{\s{\kappa^2+\sin^2\phi}}.
\end{aligned}
\ee
If $\phi\in [0,\f{\pi}{4}]$, $\sin\phi\thickapprox \phi$ and if $\phi\in [\f{\pi}{4},\f{\pi}{2}]$, $\sin\phi\thickapprox 1$, then \eqref{10.19} implies that
\bea
&&\int^{2\pi}_0|{G}(x,y)|d\phi\thickapprox \f{1}{r}\Big (\int^{\f{\pi}{4}}_0\f{d\phi}{\s{\kappa^2+\phi^2}}+\int^{\f{\pi}{2}}_{\f{\pi}{4}}\f{d\phi}{\s{\kappa^2+1}}\Big)\nn\\
&&\q\qq\qq\qq \thickapprox \f{1}{r}\Big (\int^{\f{\pi}{4\kappa}}_0\f{d\phi}{1+\phi}+1\Big)\nn\\
&&\q\qq\qq\qq\ls \f{1}{r}\ln\Big(2+\f{1}{\kappa}\Big)\nn\\
&&\q\qq\qq\qq\ls \f{1}{r}\ln\Big(2+\f{r}{|\rho-r|}\Big)\nn\\
\eea
When $1\leq|\rho-r|<r/4$, from Lemma \ref{PG}, it is easy to see that
\be
\begin{aligned}
&\q \int^{2\pi}_0|{G}(x,y)|d\phi\\
&\ls\int^{2\pi}_0\f{1}{|x'-y'|^2}d\phi\\
&\ls\int^{\pi}_0 \f{d\phi}{(\rho-r)^2+4\rho r\sin^2\phi}\\
&\ls\int^{\f{\pi}{2}}_0 \f{d\phi}{(\rho-r)^2+4\rho r\sin^2\phi}\\
& \thickapprox \f{1}{r^2}\int^{\f{\pi}{2}}_0 \f{d\phi}{\kappa^2+\sin^2\phi} \q \text{with}\q \kappa^2=\f{|\rho-r|^2}{4\rho r}\ls 1\\
& \thickapprox \f{1}{r^2}\Big(\int^{\f{\pi}{4}}_0+\int^{\f{\pi}{2}}_{\f{\pi}{4}}\Big) \f{d\phi}{\kappa^2+\sin^2\phi}.
\end{aligned}
\ee
If $\phi\in [0,\f{\pi}{4}]$, $\sin\phi\thickapprox \phi$ and if $\phi\in [\f{\pi}{4},\f{\pi}{2}]$, $\sin\phi\thickapprox 1$, then the above inequality implies that
\bea\label{2.8'}
&&\int^{2\pi}_0|{G}(x,y)|d\phi\thickapprox \f{1}{r^2}\Big (\int^{\f{\pi}{4}}_0\f{d\phi}{{\kappa^2+\phi^2}}+\int^{\f{\pi}{2}}_{\f{\pi}{4}}\f{d\phi}{{\kappa^2+1}}\Big)\nn\\
&&\q\qq\qq\qq \thickapprox \f{1}{r^2}\Big (\f{1}{\kappa}\int^{\f{\pi}{4\kappa}}_0\f{d\phi}{1+\phi^2}+1\Big)\nn\\
&&\q\qq\qq\qq\ls \f{1}{r^2}(\f{1}{\kappa}+1)\nn\\
&&\q\qq\qq\qq\ls \f{1}{r|\rho-r|}.
\eea

When $|\rho-r|<1$, from Lemma \ref{PG}, it is easy to see that
\be\label{2.9}
\begin{aligned}
&\q \int^{2\pi}_0|\p_{\rho,r}{G}(x,y)|d\phi\\
&\ls\int^{2\pi}_0\f{1}{|x'-y'|^2+|x_3-y_3|^2}d\phi\\
&\ls\int^{\pi}_0 \f{d\phi}{(\rho-r)^2+|z-l|^2+4\rho r\sin^2\phi}\\
&\ls\int^{\pi/2}_0 \f{d\phi}{(\rho-r)^2+|z-l|^2+4\rho r\sin^2\phi}\\
& \thickapprox \f{1}{r^2}\int^{\f{\pi}{2}}_0 \f{d\phi}{{\kappa^2+\sin^2\phi}} \q \text{with}\q \kappa^2=\f{|\rho-r|^2+|z-l|^2}{4\rho r}\ls 1\\
& \thickapprox \f{1}{r^2}\Big(\int^{\f{\pi}{4}}_0+\int^{\f{\pi}{2}}_{\f{\pi}{4}}\Big) \f{d\phi}{{\kappa^2+\sin^2\phi}}.
\end{aligned}
\ee
If $\phi\in [0,\f{\pi}{4}]$, $\sin\phi\thickapprox \phi$ and if $\phi\in [\f{\pi}{4},\f{\pi}{2}]$, $\sin\phi\thickapprox 1$, then \eqref{2.9} implies that
\bea
&&\int^{2\pi}_0|\p_{\rho,r}{G}(x,y)|d\phi\thickapprox \f{1}{r^2}\Big (\int^{\f{\pi}{4}}_0\f{d\phi}{{\kappa^2+\phi^2}}+\int^{\f{\pi}{2}}_{\f{\pi}{4}}\f{d\phi}{{\kappa^2+1}}\Big)\nn\\
&&\q\qq\qq\qq \thickapprox \f{1}{r^2}\Big (\f{1}{\kappa}\int^{\f{\pi}{4\kappa}}_0\f{d\phi}{1+\phi^2}+1\Big)\nn\\
&&\q\qq\qq\qq\ls \f{1}{r^2}(\f{1}{\kappa}+1)\nn\\
&&\q\qq\qq\qq\ls \f{1}{r(|\rho-r|+|z-l|)}.
\eea
When $1\leq|\rho-r|<r/4$, from Lemma \ref{PG}, it is easy to see that
\be
\begin{aligned}
&\q \int^{2\pi}_0|\p_{\rho,r}{G}(x,y)|d\phi\\
&\ls\int^{2\pi}_0\f{1}{|x'-y'|^3}d\phi\\
&\ls\int^{\pi}_0 \f{d\phi}{\big[(\rho-r)^2+4\rho r\sin^2\phi\big]^{3/2}}\\
&\ls\int^{\f{\pi}{2}}_0 \f{d\phi}{\big[(\rho-r)^2+4\rho r\sin^2\phi\big]^{3/2}}\\
& \thickapprox \f{1}{r^3}\int^{\f{\pi}{2}}_0 \f{d\phi}{\big[\kappa^2+\sin^2\phi\big]^{3/2}} \q \text{with}\q \kappa^2=\f{|\rho-r|^2}{4\rho r}\ls 1\\
& \thickapprox \f{1}{r^3}\Big(\int^{\f{\pi}{4}}_0+\int^{\f{\pi}{2}}_{\f{\pi}{4}}\Big) \f{d\phi}{\big[\kappa^2+\sin^2\phi\big]^{3/2}}.
\end{aligned}
\ee
If $\phi\in [0,\f{\pi}{4}]$, $\sin\phi\thickapprox \phi$ and if $\phi\in [\f{\pi}{4},\f{\pi}{2}]$, $\sin\phi\thickapprox 1$, then the preceding inequality implies that
\bea\label{2.12}
&&\int^{2\pi}_0|\p_{\rho,r}{G}(x,y)|d\phi\thickapprox \f{1}{r^3}\Big (\int^{\f{\pi}{4}}_0\f{d\phi}{\big[\kappa^2+\phi^2\big]^{3/2}}+\int^{\f{\pi}{2}}_{\f{\pi}{4}}\f{d\phi}{\big[\kappa^2+1\big]^{3/2}}\Big)\nn\\
&&\q\qq\qq\qq \thickapprox \f{1}{r^3}\Big (\f{1}{\kappa}\int^{\f{\pi}{4\kappa}}_0\f{d\phi}{\big[1+\phi^2\big]^{3/2}}+1\Big)\nn\\
&&\q\qq\qq\qq\ls \f{1}{r^3}(\f{1}{\kappa^2}+1)\nn\\
&&\q\qq\qq\qq\ls \f{1}{r|\rho-r|^2}.
\eea

\ef

\subsection{Decay of the velocity $u$:  $|u| \ls   \big(\f{\ln r}{r}\big)^{1/2}.$}

In this subsection, we prove the estimate.
\be
\lab{u-1/2}
|u|\ls \big(\f{\ln r}{r}\big)^{1/2}.
\ee

The proof is based on the following $Brezis-Gallouet$ inequality \cite{BG:1} and its refinement. together with scaling techniques.

\begin{lemma}\label{l3.1}
Let $f\in H^2(\O)$ where $\O \subset \bR^2$. Then there exists a constant $C_{\O}$, depending only on $\O$, such that
\be
\lab{BGinq}
\|f\|_{L^\i(\O)}\leq C_{\O}\|f\|_{H^1(\O)}\log^{1/2} \big(e+\f{\|\Dl f\|_{L^2(\O)}}{\|f\|_{H^1(\O)}}\big).
\ee
\end{lemma}
\begin{remark}
Usually we will  use the following inequality instead of \eqref{BGinq}
\be\label{e2.15e}
\|f\|_{L^\i(\O)}\leq C_{\O}(1+\|f\|_{H^1(\O)})\log^{1/2} \big(e+\|\Dl f\|_{L^2(\O)}\big).
\ee
\end{remark}

Using Lemma \ref{l3.1}, we can prove the following refined $Brezis-Gallouet$ inequality whose constant is independent of the thinness of the domain.
A price to pay is that the functions need to have zero boundary value in the $z$ direction or
mean zero in the $z$ direction.  This property will be needed later.

\begin{lemma}
\lab{leRBG}
Set
\be
\bar{\mathcal{D}}_0=\{(r,z): 1/2<r<2,\ z\in[0,1/\la]\}.\nn
\ee
Then if $f \in H^2(\bar{\mathcal{D}}_0)$ satisfies
\be\label{e2.19}
\int^{1/\la}_0 fdz=0\q \text{or} \q f|_{z=0,1/\la}=0,
\ee
we have
\be\label{2.19}
\|f\|_{L^\i(\bar{\mathcal{D}}_0)}\leq C_0(1+\|\na f\|_{L^2(\bar{\mathcal{D}}_0)})\log^{1/2}(e+\f{1}{\la}\|\Dl f\|_{L^2(\bar{\mathcal{D}}_0)}),
\ee
where $C_0$ is independent of $\la$.
\end{lemma}
\pf
 Note that we can not simply make zero extension for $f$ outside of the domain and apply the
 regular Brezis-Gallouet inequality. The reason is that the extended function may not be in $H^2$.

Let $\t{f}(\t{r},\t{z})=f(\t{r}/\la,\t{z}/\la)$ where $(\t{r},\t{z})\in \bar{\mathcal{D}}_{0,\la}$ and
\be
\bar{\mathcal{D}}_{0,\la}=\{(r,z):\la/2<r<2 \la,\ z\in[0,1]\}.\nn
\ee
Using \eqref{e2.15e}, we get
\bea\label{2.20}
&&\|f(r,z)\|_{L^\i(\bar{\mathcal{D}}_0)}\nn\\
&=&\|\t{f}(\t{r},\t{z})\|_{L^\i(\bar{\mathcal{D}}_{0,\la})}\nn\\
&\leq& C_0(1+\|\t{f}\|_{H^1(\bar{\mathcal{D}}_{0,\la})})\log^{1/2} \big(e+\|\Dl \t{f}\|_{L^2(\bar{\mathcal{D}}_{0,\la})}\big)\nn\\
&=& C_0(1+\|\na\t{f}\|_{L^2(\bar{\mathcal{D}}_{0,\la})}+\|\t{f}\|_{L^2(\bar{\mathcal{D}}_{0,\la})})\log^{1/2} \big(e+\|\Dl \t{f}\|_{L^2(\bar{\mathcal{D}}_{0,\la})}\big).
\eea We mention that the constant $C_0$ is independent of $\la$. The reason is that we can first extend the function $\t{f}$ to be a $H^2$ function in the whole $(r, z)$ space.  From the proof of the original
Brezis-Gallouet inequality, we know the constant relies only on the $H^2$ extension property of functions in a domain. The extension property only depends on the thickness of the original rectangle, which is scaled to $1$. Alternatively, one can also just pick a point in
$\bar{\mathcal{D}}_{0,\la}$ and apply the usual Brezis- Gallouet inequality in a unit ball centered at this point.

By the variable change and relationship between $f$ and $\t{f}$. we can get
\be\label{2.21}
\|\na\t{f}\|_{L^2(\bar{\mathcal{D}}_{0,\la})}=\|\na{f}\|_{L^2(\bar{\mathcal{D}}_{0})}, \|\t{f}\|_{L^2(\bar{\mathcal{D}}_{0,\la})}=\la\|{f}\|_{L^2(\bar{\mathcal{D}}_{0})},
\|\Dl \t{f}\|_{L^2(\bar{\mathcal{D}}_{0,\la})}=\f{1}{\la}\|\Dl f \|_{L^2(\bar{\mathcal{D}}_{0})}.
\ee
Inserting \eqref{2.21} into \eqref{2.20}, we get
\bea\label{2.22}
&&\|f(r,z)\|_{L^\i(\bar{\mathcal{D}}_0)}\nn\\
&\ls& C_0(1+\|\na{f}\|_{L^2(\bar{\mathcal{D}}_{0})}+\la\|{f}\|_{L^2(\bar{\mathcal{D}}_{0})})\log^{1/2} \big(e+\f{1}{\la}\|\Dl f \|_{L^2(\bar{\mathcal{D}}_{0})}\big).
\eea

Now if $f$ satisfies \eqref{e2.19}, by using Poincar\'e inequality, we have
\be
\la\|{f}\|_{L^2(\bar{\mathcal{D}}_{0})}\leq C\|\na{f}\|_{L^2(\bar{\mathcal{D}}_{0})},\nn
\ee
where $C$ is independent of $\la$. At last, combination of the above inequality and \eqref{2.22} indicates \eqref{2.19}.\ef

Now we can prove (\ref{u-1/2}), the decay of $u$.

Fixing $x_0 \in \bR^2 \times [0, 1]$ such that $|x'_0|=  r_0$ is large. Without loss of generality, we can assume, in the cylindrical coordinates, that $x_0=(r_0, 0, 0)$, i.e. $z_0=0$, $\theta_0=0$.
Consider the scaled  solution
\[
\t {u} ( \t{x}) = r_0 u( r_0 \t{x})
\]which is also axially symmetric. Hence $\t {u}$ can be regarded as a two variable function  of the scaled variables $\t{r}, \t{z}$. Consider the two dimensional domain
\[
\t D =\{ (\t{r}, \t{z}) \, |,
1/2 \le \t{r} \le 2, \,  |\t{z}| \le 1/r_0 \}.
\]Then for $\t{u}=\t{u}(\t{r}, \t{z})$, we have $\t{u}(1, 0)=r_0 u(x_0)$.

Recall that $u$ satisfies the Dirichlet boundary condition.
Applying the refined Brezis-Gallouet inequality (Lemma \ref{leRBG})
 on $\t{D}$, after a simple adjustment on constants, we can find an absolute constant $C$ such that
\[
\al
| \t{u}(1, 0) | &\le C \left[ \left( \int_{\t{D}} |\t{\nabla} \t{u} |^2 d\t{r}d\t{z} \right)^{1/2}
 +1  \right] \times \\
&\qquad \log^{1/2}\left[ \left( \int_{\t{D}} |\t{\Delta} \t{u} |^2 d\t{r}d\t{z} \right)^{1/2} +
\left(\int_{\t{D}} | \t{u} |^2 d\t{r}d\t{z}\right)^{1/2} + e  \right],
\eal
\]where $\t \nabla = (\partial_{\t{r}}, \partial_{\t{z}})$ and
$\t \Delta = \partial^2_{\t{r}} + \partial^2_{\t{z}}$.  By  the assumption that
$1/2 \le \t{r} \le 2$, we see that
\[
\al
| \t{u}(1, 0) | &\le C \left[ \left( \int_{\t{D}} |\t{\nabla} \t{u} |^2 \t{r} d\t{r}d\t{z} \right)^{1/2}  +1  \right] \times \\
&\qquad \log^{1/2}\left[ \left( \int_{\t{D}} |\t{\Delta} \t{u} |^2 \t{r} d\t{r}d\t{z} \right)^{1/2} +
\left(\int_{\t{D}} | \t{u} |^2 \t{r} d\t{r}d\t{z}\right)^{1/2} + e  \right].
\eal
\]Now we can scale this inequality back to the original solution $u$ and variables $r=r_0 \t{r}$ and $z=r_0 \t{z}$ to get
\[
\al
r_0 |u(x_0)| &\le C  \left[ \sqrt{r_0} \left( \int_{D_0} |\nabla u |^2 r drdz \right)^{1/2}  +1  \right] \times \\
&\qquad \log^{1/2}\left[ r^{3/2}_0 \left( \int_{D_0} (|\p^2_r u|^2 + |\p^2_z u|^2) r drdz \right)^{1/2} + r^{-1/2}_0
\left(\int_{D_0} | u |^2 r drdz \right)^{1/2} + e  \right],
\eal
\]where
\[
 D_0 =\{ (r, z) \, |
r_0/2 \le r \le 2 r_0, \,  0 \le z  \le 1 \}.
\]By condition (\ref{dr2r}), this proves the claimed decay of velocity. Note that by our assumption in the theorem, the solution $u$ is globally bounded and then  it's not hard to prove that the first and second derivatives of $u$ are also bounded.

\subsection{Decay of $\boldsymbol{w^r,\ w^z:\ |(w^r,w^z)|\ls r^{-1}\ln r}$.}

\q First we investigate the boundary conditions for $w^r, w^z$. We will show that
\be\label{2.12'}
\p_z w^r\big|_{z=0,1}=0,\ w^z\big|_{z=0,1}=0.
\ee These $0$ boundary values allow us to work on the vorticity equation with ease.

From \eqref{1.4} and \eqref{e1.4}, we see
\[
\p_z w^r=-\p^2_z u^\th=-(u^r\p_r+u^z\p_z)u^\th-\f{u^ru^\th}{r}+(\p^2_r+\f{1}{r}\p_r-\f{1}{r^2})u^\th.
\]The Dirichlet boundary condition on $u^r,u^\th, u^z$ indicates that
\be
\p_z w^r\big|_{z=0,1}=-(u^r\p_r+u^z\p_z)u^\th-\f{u^ru^\th}{r}+
\left(\p^2_r+\f{1}{r}\p_r-\f{1}{r^2} \right)u^\th\big|_{z=0,1}=0. \nn
\ee

Also from \eqref{e1.4}, we have
\be
w^z\big|_{z=0,1}=\f{1}{r}\p_r \left( ru^\th \right)\big|_{z=0,1}=0. \nn
\ee

So we can do integration by parts for the first and third equation of \eqref{1.5} without any boundary terms coming out. Besides we have
\be\label{2.13}
\int^1_0 w^rdz=\int^1_0 -\p_z u^\th dz=-u^\th(r,1)+u^\th(r,0)=0.
\ee

We pick a point $x_0 \in \bR^3$ such that $|x'_0| \equiv \la $ is large and carry out the
scaling for the velocity and vorticity:
\bea
&&\t{u}(\t{x})=\la u(\la \t{x})=\la u(x)\nn\\
&&\t{w}(\t{x})=\la^2w(\la \t{x})=\la^2w(x)\nn
\eea
where $\t{x}=\f{x}{\la}$.

For simplification of notation, we temporarily drop the $``\sim"$ symbol when computations take place under the scaled sense. Define the domains
\be
\mathcal{D}_1=\{(r,\th,z):\f{1}{2}<r<\f{3}{2},\ 0\leq \th\leq 2\pi,\ z\in[0,1/\la]\}\nn
\ee
and
\be
\mathcal{D}_2=\{(r,\th,z):\f{3}{4}<r<\f{5}{4},\ 0\leq \th\leq 2\pi,\ z\in[0,1/\la]\}.\nn
\ee

Let $\psi=\psi(y')$ be a cut-off function satisfying $\supp \psi(y')\subset \mathcal{D}_1$ and $\psi(y')=1$ for $y\in\mathcal{D}_2$ such that the gradient of $\psi$ is bounded. Note that $\psi$ is independent of the $z$ variable. Now testing the first and third one of the vorticity equations \eqref{1.5} with $w^r\psi^2$ and $w^z\psi^2$ respectively, we have

\bea
&&\q-\int_{\mathcal{D}_1} w^r\psi^2(\Dl-\f{1}{r^2})w^r dy\nn\\
&&=-\int_{\mathcal{D}_1}\left[ (u^r\p_r+u^z\p_z)w^r\cdot w^r\psi^2- (w^r\p_r+w^z\p_z)u^r\cdot w^r\psi^2 \right]dy. \nn
\eea
\bea
&&-\int_{\mathcal{D}_1} w^z\psi^2\Dl w^z dy\nn\\
&&=-\int_{\mathcal{D}_1}\left[(u^r\p_r+u^z\p_z)w^z\cdot w^z\psi^2-(w^r\p_r+w^z\p_z)u^z\cdot w^z\psi^2 \right]dy. \nn
\eea

Then we have
\be\label{4.1}
\begin{aligned}
&\q\ \int_{\mathcal{D}_1} \Big(|\nabla(w^r\psi)|^2+\f{(w^r)^2\psi^2}{r^2}\Big) dy  \\
&=\int_{\mathcal{D}_1}\Big((w^r)^2|\nabla\psi|^2-\f{1}{2}\psi^2(u^r\p_r+u^z\p_z)(w^r)^2\\
&\qq\qq\qq\qq\qq+(w^r)^2\psi^2\p_ru^r+w^rw^z\psi^2\p_zu^r\Big )dy\\
&=\int_{\mathcal{D}_1}\Big((w^r)^2|\nabla\psi|^2+\f{1}{2}(w^r)^2(u^r\p_r+u^z\p_z)\psi^2\\
&\qq\qq\qq\qq\qq-2u^rw^r\psi\p_r(w^r\psi)-(w^r \psi)^2 \frac{u^r}{r}-u^r\p_z(w^r\psi w^z\psi)\Big )dy\\
&\leq C(1+\|(u^r,u^z)\|_{L^\i(\mathcal{D}_1)})\|w^r\|^2_{L^2(\mathcal{D}_1)}+\f{1}{4}\big(\|\nabla(w^r\psi)\|^2_{L^2(\mathcal{D}_1)}+\|\nabla(w^z\psi)\|^2_{L^2(\mathcal{D}_1)}\big)\\
&\qq\qq\qq\qq\qq+C\|u^r\|^2_{L^\i(\mathcal{D}_1)}\big(\|w^r\|^2_{L^2(\mathcal{D}_1)}+\|w^z\|^2_{L^2(\mathcal{D}_1)}\big)\\
&\leq C(1+\|(u^r,u^z)\|^2_{L^\i(\mathcal{D}_1)})\|(w^r,w^z)\|^2_{L^2(\mathcal{D}_1)}+\f{1}{4}\|\big(\nabla(w^r\psi),\nabla(w^z\psi)\big)\|^2_{L^2(\mathcal{D}_1)},
\end{aligned}
\ee
and
\be\label{4.3}
\begin{aligned}
&\q\ \int_{\mathcal{D}_1} |\nabla(w^z\psi)|^2 dy  \\
&=\int_{\mathcal{D}_1}\Big((w^z)^2|\nabla\psi|^2-\f{1}{2}\psi^2(u^r\p_r+u^z\p_z)(w^z)^2\\
&\qq\qq\qq\qq\qq+w^rw^z\psi^2\p_ru^z+(w^z\psi)^2\p_zu^z\Big )dy\\
&=\int_{\mathcal{D}_1}\Big((w^z)^2|\nabla\psi|^2+\f{1}{2}(w^z)^2(u^r\p_r+u^z\p_z)\psi^2\\
&\qq\qq\qq\qq\qq-2u^zw^z\psi\p_z(w^z\psi)-u^z\p_r(w^r\psi w^z\psi) -w^z w^r \psi^2 \frac{u^z}{r} \Big )dy\\
&\leq C(1+\|(u^r,u^z)\|_{L^\i(\mathcal{D}_1)})\|w^z\|^2_{L^2(\mathcal{D}_1)}+\f{1}{4}\big(\|\nabla(w^r\psi)\|^2_{L^2(\mathcal{D}_1)}+
\|\nabla(w^z\psi)\|^2_{L^2(\mathcal{D}_1)}\big)\\
&\qq\qq\qq\qq\qq+C\|(u^r, u^z)\|^2_{L^\i(\mathcal{D}_1)}\big(\|w^r\|^2_{L^2(\mathcal{D}_1)}+\|w^z\|^2_{L^2(\mathcal{D}_1)}\big)\\
&\leq C(1+\|(u^r,u^z)\|^2_{L^\i(\mathcal{D}_1)})\|(w^r,w^z)\|^2_{L^2(\mathcal{D}_1)}+
\f{1}{4}\|\big(\nabla(w^r\psi),\nabla(w^z\psi)\big)\|^2_{L^2(\mathcal{D}_1)}.
\end{aligned}
\ee
From \eqref{4.1} and \eqref{4.3} we obtain
\be \label{4.41}
\begin{aligned}
&\q\ \|\big(\nabla w^r, \nabla w^z\big)\|^2_{L^2(\mathcal{D}_2)}\\
 & \leq  C(1+\|(u^r,u^z)\|_{L^\i(\mathcal{D}_1)})^2\|(w^r,w^z)\|^2_{L^2(\mathcal{D}_1)}.
 \end{aligned}
\ee

Since our scaled $w^r, w^z$ satisfies \eqref{e2.19}, we apply \eqref{2.19} to get
\bea\label{2.24}
&&\|(w^r,w^z)\|_{L^\i(\bar{\mathcal{D}}_0)}\nn\\
&\leq& C_0(1+\|\na (w^r,w^z)\|_{L^2(\bar{\mathcal{D}}_0)})\log^{1/2}(e+\f{1}{\la}\|\Dl (w^r,w^z)\|_{L^2(\bar{\mathcal{D}}_0)}).
\eea
Now inserting \eqref{4.41} into \eqref{2.24} implies
\bea\label{2.25}
&&\|(w^r,w^z)\|_{L^\i(\bar{\mathcal{D}}_0)}\nn\\
&\leq& C_0\Big(1+(1+\|(u^r,u^z)\|_{L^\i(\mathcal{D}_1)})\|(w^r,w^z)\|_{L^2(\mathcal{D}_1)}\Big)\times\nn\\
&&\log^{1/2}(e+\f{1}{\la}\|\Dl (w^r,w^z)\|_{L^2(\bar{\mathcal{D}}_0)}).
\eea
Now scaling back, to the domains
\be
\mathcal{D}_{0, \la}=\{(r,\th,z):\la-1<r<\la+1,\ 0\leq \th\leq 2\pi,\ z\in[0,1]\},\nn
\ee
\be
\mathcal{D}_{1, \la}=\{(r,\th,z):\f{\la}{2}<r<\f{3\la}{2},\ 0\leq \th\leq 2\pi,\ z\in[0,1]\},\nn
\ee
we can get
\bea\label{2.26}
&&\la^2\|(w^r,w^z)\|_{L^\i({\mathcal{D}}_{0,\la})}\nn\\
&\leq& C_0\Big(1+(1+\la\|(u^r,u^z)\|_{L^\i(\mathcal{D}_{1,\la})})\la^{1/2}\|(w^r,w^z)\|_{L^2(\mathcal{D}_{1,\la})}\Big)\times\nn\\
&&\log^{1/2}(e+\la^{3/2}\|\Dl (w^r,w^z)\|_{L^2(\bar{\mathcal{D}}_{0,\la})}).
\eea

From (\ref{u-1/2}) we have $|(u^r,u^z)(x)|\ls \f{(\ln r)^{1/2}}{r^{1/2}}$ with $|x'|=r$, therefore \eqref{2.26} implies that
\be\label{2.27}
|(w^r,w^z)(x_0)|\ls \la^{-1}\ln \la,
\ee with $\la=|x'_0|$.

\subsection{Decay of $\boldsymbol{\p_z w^r}$ and $\boldsymbol{\p_r w^z}$: $\boldsymbol{|(\p_z w^r,\p_r w^z)|\ls r^{-3/2}(\ln r)^{3/2}}$.}

Now we use the scaling technique and refined $Brezis-Gallouet$ inequality to prove the decay of $J:=\p_z w^r$ and $\O=\p_r w^z$. First we observe that the boundary conditions of $J,\O$ are $0$:
\be\label{2.28}
J|_{z=0,1}=\O|_{z=0,1}=0.
\ee These follow from the relation that $w^r=-\p_z u^\th$, and $w^z= \p_r u^\th + u^\th/r$ and
the equation for $u^\th$.
From \eqref{1.5}, we see $(J,\O)$ satisfy
\be\label{2.29}
\lt\{
\begin{aligned}
&(\Dl-\f{1}{r^2})J=\p_z\big[(u^r\p_r+u^z\p_z)w^r\big]-\p_z\big[(w^r\p_r+w^z\p_z)u^r\big],\\
&(\Dl-\f{1}{r^2})\O=\p_r\big[(u^r\p_r+u^z\p_z)w^z\big]-\p_r\big[(w^r\p_r+w^z\p_z)u^z\big].
\end{aligned}
\rt.
\ee Here we are still working on scaled functions and domains without using the tilde notation, unless stated otherwise.

Set
\be
\mathcal{D}_3=\{(r,\th,z):\f{7}{8}<r<\f{9}{8},\ 0\leq \th\leq 2\pi,\ z\in[0,1/\la]\}.\nn
\ee
Let $\psi(y')$ be a cut-off function satisfying $\supp \psi(y')\subset \mathcal{D}_2$ and $\psi(y')=1$ for $y\in\mathcal{D}_3$ such that the gradient of $\psi$ is bounded. Now testing \eqref{2.28} with $J\psi^2$ and $\O\psi^2$ respectively, we have

\bea
&&\q-\int_{\mathcal{D}_2} J\psi^2(\Dl-\f{1}{r^2})J dy\nn\\
&&=\int_{\mathcal{D}_2}\Big(-\p_z\big[(u^r\p_r+u^z\p_z)w^r\big]+\p_z\big[(w^r\p_r+w^z\p_z)u^r\big] \Big)J\psi^2dy. \nn
\eea
\bea
&&-\int_{\mathcal{D}_2} \O\psi^2(\Dl-\f{1}{r^2})\O dy\nn\\
&&=\int_{\mathcal{D}_2}\Big(-\p_r\big[(u^r\p_r+u^z\p_z)w^z\big]+\p_r\big[(w^r\p_r+w^z\p_z)u^z\big] \Big)\O\psi^2dy. \nn
\eea

Then we have
\be
\begin{aligned}
&\q\ \int_{\mathcal{D}_2} \Big(|\nabla(J\psi)|^2+\f{J^2\psi^2}{r^2}\Big) dy  \\
&=\int_{\mathcal{D}_2}\Big(J^2|\nabla\psi|^2+(u^r\p_r+u^z\p_z)w^r \p_z(J\psi^2)-(w^r\p_r+w^z\p_z)u^r \p_z(J\psi^2)\Big )dy\\
&\ls \|J\|^2_{L^2(\mathcal{D}_2)}+\|(u^r,u^z)\|^2_{L^\i(\mathcal{D}_2)}\|\na w^r\|^2_{L^2(\mathcal{D}_2)}\nn\\
&\qq+\|(w^r,w^z)\|^2_{L^\i(\mathcal{D}_2)}\|\na u^r\|^2_{L^2(\mathcal{D}_2)}+\f{1}{2}\|\nabla(J\psi)\|^2_{L^2(\mathcal{D}_2)},\nn
\end{aligned}
\ee
and
\be
\begin{aligned}
&\q\ \int_{\mathcal{D}_2} \Big(|\nabla(\O\psi)|^2+\f{\O^2\psi^2}{r^2}\Big) dy  \\
&=\int_{\mathcal{D}_2}\Big(\O^2|\nabla\psi|^2+(u^r\p_r+u^z\p_z)w^z \p_r(\O\psi^2)-(w^r\p_r+w^z\p_z)u^z \p_r(\O\psi^2)\Big )dy\\
&=\int_{\mathcal{D}_2}\Big(\O^2|\nabla\psi|^2+(u^r\p_r+u^z\p_z)w^z (\p_r(\O\psi)\psi+\O\psi\p_r\psi)\nn\\
&\qq\qq\qq-(w^r\p_r+w^z\p_z)u^z \big(\p_r(\O\psi)\psi+\O\psi\p_r\psi\big)\Big )dy\\
&\ls \|\O\|^2_{L^2(\mathcal{D}_2)}+\|(u^r,u^z)\|^2_{L^\i(\mathcal{D}_2)}\|\na w^z\|^2_{L^2(\mathcal{D}_2)}\nn\\
&\qq+\|(w^r,w^z)\|^2_{L^\i(\mathcal{D}_2)}\|\na u^z\|^2_{L^2(\mathcal{D}_2)}+\f{1}{2}\|\nabla(\O\psi)\|^2_{L^2(\mathcal{D}_2)}.\nn
\end{aligned}
\ee
The above two inequalities indicate that
\be\label{2.30}
\begin{aligned}
&\q\ \|\big(\nabla J, \nabla \O\big)\|^2_{L^2(\mathcal{D}_3)}\\
 & \ls  (1+\|(u^r,u^z)\|^2_{L^\i(\mathcal{D}_2)})\|(\na w^r,\na w^z)\|^2_{L^2(\mathcal{D}_2)}\\
 &\q+\|(w^r,w^z)\|^2_{L^\i(\mathcal{D}_2)}\|\na (u^r,u^z)\|^2_{L^2(\mathcal{D}_2)}.
 \end{aligned}
\ee

Inserting \eqref{4.41} into \eqref{2.30} implies
\be\label{2.31}
\begin{aligned}
&\q\ \|\big(\nabla J, \nabla \O\big)\|^2_{L^2(\mathcal{D}_3)}\\
 & \ls  (1+\|(u^r,u^z)\|_{L^\i(\mathcal{D}_1)})^4\|(w^r,w^z)\|^2_{L^2(\mathcal{D}_1)}\\
 &\q+\|(w^r,w^z)\|^2_{L^\i(\mathcal{D}_1)}\|\na (u^r,u^z)\|^2_{L^2(\mathcal{D}_1)}.
 \end{aligned}
\ee
We apply \eqref{2.19} to get
\bea
&&\|(J,\O)\|_{L^\i(\bar{\mathcal{D}}_0)}\nn\\
&\leq& C_0(1+\|\na (J,\O)\|_{L^2(\bar{\mathcal{D}}_0)})\log^{1/2}(e+\f{1}{\la}\|\Dl (J,\O)\|_{L^2(\bar{\mathcal{D}}_0)}),\nn\\
&\leq& C_0\Big[1+(1+\|(u^r,u^z)\|_{L^\i(\mathcal{D}_1)})^2\|(w^r,w^z)\|_{L^2(\mathcal{D}_1)}\nn\\
&&+\|(w^r,w^z)\|_{L^\i(\mathcal{D}_1)}\|\na (u^r,u^z)\|_{L^2(\mathcal{D}_1)}\Big]\log^{1/2}(e+\f{1}{\la}\|\Dl (J,\O)\|_{L^2(\bar{\mathcal{D}}_0)}).
\eea
Now scaling back to the original functions and domains $\mathcal{D}_{0,\la},\ \mathcal{D}_{1,\la}$, we deduce
\bea
&&\la^3\|(J,\O)\|_{L^\i(\mathcal{D}_{0,\la})}\nn\\
&\leq& C_0\Big[1+(1+\la \|(u^r,u^z)\|_{L^\i(\mathcal{D}_{1,\la})})^2\la^{1/2}\|(w^r,w^z)\|_{L^2(\mathcal{D}_{1,\la})}\nn\\
&&+\la^2\|(w^r,w^z)\|_{L^\i(\mathcal{D}_{1,\la})}\la^{1/2}\|\na (u^r,u^z)\|_{L^2(\mathcal{D}_{1,\la})}\Big]\log^{1/2}(e+\la^{5/2}\|\Dl (J,\O)\|_{L^2(\mathcal{D}_{0,\la})}). \nn
\eea

The above inequality implies that
\be\label{2.33'}
\|(J,\O)(x_0)\|_{L^\i(\mathcal{D}_{0,\la})}\ls \la^{-3/2}(\ln\la)^{3/2}, \quad |x_0'| =\la.
\ee This is the claimed decay for $\p_z w^r$ and $\p_r w^z$.

\subsection{Decay of $\boldsymbol{u^\th:|u^\th|\ls r^{-3/2}(\ln r)^{5/2}}$ and vanishing.}

In this subsection we prove that $u^\theta$ decays faster than first order and then use the maximum principle to conclude that $u^\th=0$.

\q Using the $Biot-Savart$ law, for a cut off functions $\psi=\psi(x')$, which is independent of $x_3$, we know, for any smooth, divergence free vector field $v$, that
\be
-\Dl(v \psi)=\psi\na\times (\na \times v) -2\na\psi\cdot\na v -v\Dl\psi.\nn
\ee

Then using the Green's function on $\bR^2\times [0,1]$ with Dirichlet boundary, we have
\be\label{e2.35}
\begin{aligned}
\psi v(x)&=\int^1_{0}\int_{\bR^2}G(x,y)\psi\na\times (\na \times v) dy\\
&-2\int^1_{0}\int_{\bR^2}G(x,y)\na\psi\cdot\na v dy-\int^1_0\int_{\bR^2}G(x,y)(\Dl\psi) v dy.
\end{aligned}
\ee

If we write $v=u^\th e_\th$, then $\na\times v=w^re_r+w^ze_z$.  Let $x=(r\cos\th,r\sin\th,z)$, $y=(\rho\cos\phi,\rho\sin\phi,l)$. Then from \eqref{e2.35}, we have
\bea
&&\psi u^\th(x)=\int^1_{0}\int_{\bR^2}G(x,y)\psi(\na\times (w^\rho e_\rho+w^le_l))\cdot e_\th dy\nn\\
&&\qq-2\int^1_{0}\int_{\bR^2}G(x,y)(\na\psi\cdot\na v)\cdot e_\th dy-\int^1_{0}\int_{\bR^2}G(x,y)(\Dl\psi) v\cdot e_\th d y\nn\\
&&\q\ =\int^{1}_{0}\int^{2\pi}_0\int^\i_0G(x,y)\psi (\p_l w^\rho-\p_\rho w^l) \cos(\phi-\th)\rho d\rho d\phi dl\nn\\
&&\qq -2\int^{1}_{0}\int^{2\pi}_0\int^\i_0G(x,y)\p_\rho\psi\p_\rho u^\phi \cos(\phi-\th) \rho d\rho d\phi dl\nn\\
&&\qq -\int^{1}_{0}\int^{2\pi}_0\int^\i_0G(x,y)(\p^2_\rho\psi+\f{1}{\rho}\p_\rho\psi) u^\phi \cos(\phi-\th) \rho d\rho d\phi dl.\nn
\eea
By setting $\th=0$ and integration by part for the second line of the above equality, we have
\be\label{2.34}
\begin{aligned}
\psi u^\th(x)&=\underbrace{\int^{1}_{0}\int^\i_0\Big(\int^{2\pi}_0{G}(x,y)\cos \phi d\phi\Big)(\p_lw^\rho-\p_\rho w^l) \psi\rho d\rho  dl}_{I_1}\\
&\q\ +2\underbrace{\int^{1}_{0}\int^\i_0\Big(\int^{2\pi}_0 \p_\rho{G}(x,y)\cos\phi d\phi\Big)\p_\rho\psi u^\phi  \rho d\rho dl}_{I_2}\\
&\q\ +\underbrace{\int^{1}_{0}\int^\i_0\Big(\int^{2\pi}_0 {G}(x,y)\cos\phi d\phi\Big)(\p^2_\rho\psi+\f{1}{\rho}\p_\rho\psi) u^\phi \rho d\rho d\phi dl}_{I_3}.\\
\end{aligned}
\ee Select the cut-off function $\psi$ such that its support is contained in the annulus
$B(0, 5 r/4)-B(0, 3 r/4)$ with $r = |x'|$. Then,

\bea
&&|I_1|\ls \int^1_{0}\Big(\int_{|\rho-r|<1}+\int_{1\leq |\rho-r|\leq \f{1}{4}r}\Big)\int^{2\pi}_0|{G}(x,y)|d\phi(| \p_\rho w^l|+|\p_l w^\rho|)\rho d\rho dl\nn\\
&&\qq \ls \sup\limits_{(\rho,l)\in\text{supp} \psi} (| \p_\rho w^l|+|\p_l w^\rho|)\Bigg\{\int^1_{0}\int_{|\rho-r|\leq 1} \ln(1+\f{r}{|\rho-r|})d\rho dl\nn\\
&&\qq\qq\qq\qq+\int^1_{0}\int_{1<|\rho-r|\leq \f{r}{4}} \f{1}{|\rho-r|}d\rho dl\Bigg\}\nn\\
&&\qq \ls r^{-3/2}(\ln r)^{3/2}\underbrace{\Big(\int^1_{0}\ln(1+\f{r}{s})ds+\int^{\f{r}{4}}_1\f{1}{s}ds\Big)}_{J}.\nn
\eea
Since \be
J\ls \ln r, \nn
\ee we have
\be\label{2.35}
I_1\ls r^{-3/2}(\ln r)^{5/2}.
\ee

\bea\label{2.36}
&&|I_2|\ls \int^1_{0}\Big(\int_{|\rho-r|<1}+\int_{1\leq |\rho-r|\leq \f{1}{4}r}\Big)\int^{2\pi}_0|\p_\rho{G}(x,y)|d\phi|\p_\rho \psi||u^\phi|\rho d\rho dl\nn\\
&&\qq \ls \sup\limits_{(\rho,l)\in\text{supp} \p_\rho\psi} |u^\phi|\Bigg\{\int^1_{0}\int_{|\rho-r|\leq 1} \f{1}{r(|\rho-r|+|z-l|)}d\rho dl\nn\\
&&\qq\qq\qq\qq+\int^1_{0}\int_{1<|\rho-r|\leq \f{r}{4}} \f{1}{r|\rho-r|^2}d\rho dl\Bigg\}\nn\\
&&\qq \ls r^{-1/2}(\ln r)^{1/2}\Big(\f{1}{r}\int^1_{0}\int^1_0
\f{1}{s+t}dsdt+\f{1}{r}\int^{\f{r}{4}}_1\f{1}{s^2}ds\Big)\nn\\
&&\qq \ls r^{-3/2}(\ln r)^{1/2}.
\eea

\bea\label{2.37}
&&|I_3|\ls \int^1_{0}\Big(\int_{|\rho-r|<1}+\int_{1\leq |\rho-r|\leq \f{1}{4}r}\Big)\int^{2\pi}_0|{G}(x,y)|d\phi(|\p^2_\rho \psi|+\f{|\p_\rho\psi|}{\rho})|u^\phi|\rho d\rho dl\nn\\
&&\qq \ls \f{1}{r}\sup\limits_{(\rho,l)\in\text{supp} \p_\rho\psi} |u^\phi|\Bigg\{\int^1_{0}\int_{|\rho-r|\leq 1} \f{1}{r}\ln(1+\f{r}{|\rho-r|})d\rho dl\nn\\
&&\qq\qq\qq\qq+\int^1_{0}\int_{1<|\rho-r|\leq \f{r}{4}} \f{1}{r|\rho-r|}d\rho dl\Bigg\}\nn\\
&&\qq \ls r^{-3/2}(\ln r)^{1/2}\Big(\f{1}{r}\int^1_{0}\ln(1+\f{r}{s})ds+\f{1}{r}\int^{\f{r}{4}}_1\f{1}{s}ds\Big)\nn\\
&&\qq \ls r^{-5/2}(\ln r)^{3/2}.
\eea

Combination of \eqref{2.34}$-$\eqref{2.37}, we can get
\be\label{2.38}
|u^\th|\ls r^{-3/2}(\ln r)^{5/2},
\ee
which induces the vanishing of $u^\th$. Indeed, the function $\Gamma = r  u^\th$ is known to satisfy the equation
\[
\Delta \Gamma - \frac{2}{r} \p_r \Gamma - (u^r \p_r + u^z \p_z ) \Gamma = 0.
\]We can regard this equation as a 2 dimensional one for the variables $(r, z)$. i.e.
\[
(\p^2_r + \p^2_z) \Gamma - \frac{1}{r} \p_r \Gamma - (u^r \p_r + u^z \p_z ) \Gamma = 0.
\]

From (\ref{2.38}), we see that
\[
\lim_{r \to \infty} \Gamma(r, z) =0
\]uniformly for all $z$. Also $\Gamma=0$ on the boundary of the 2 dimensional domain
\[
\{ (r, z ) \, | \, r \ge 0 , z \in [0, 1]\}.
 \]Hence, $\Gamma$ is identically $0$ and thus $u^\th \equiv 0$. Otherwise  there must be an interior maximum or minimum, violating the strong maximum principle which applies in the interior, since the coefficients of the equation are regular there. The theorem is proven. \qed

\medskip

\section*{Appendix}
\renewcommand\theequation{A.\arabic{equation}}

In the Appendix, we are devoted to proving the Lemma \ref{PG} which is equivalent to the following estimates
\be
|G(x,y)|\ls\lt\{
\begin{aligned}
&\f{1}{|x'-y'|+|x_3-y_3|}\qq  |x'-y'|<1,\\
&\f{1}{|x'-y'|^2} \qq\qq\q\q\ \ |x'-y'|>1,
\end{aligned}
\rt.\label{A.1}
\ee
and
\be
|\p_{x',y'}G(x,y)|\ls\lt\{
\begin{aligned}
&\f{1}{|x'-y'|^2+|x_3-y_3|^2} \qq  |x'-y'|<1,\\
&\f{1}{|x'-y'|^3} \qq\qq\q\q\q \ |x'-y'|>1.
\end{aligned}
\rt.\label{A.2}
\ee
The proof is by direct computation. Alternatively, one can also first consider the heat kernel with Dirichlet boundary condition on $\bR^2 \times [0, 1]$, which has fast decay in time.   Then
one can integrate out the time variable to obtain the Green's function estimate. This was the route taken in \cite{CPZ:1}.
We only show the proof of \eqref{A.1}, and the proof of \eqref{A.2} will be essentially the same and we omit the details.

For simplification of notation and without cause of confusion, we denote
\[
k_{n,-,+}=\s{|x'-y'|^2+|x_3-y_3+2n|^2},\q k_{n,+,-}=\s{|x'-y'|^2+|x_3+y_3-2n|^2},
\]
\[
k_{n,-,-}=\s{|x'-y'|^2+|x_3-y_3-2n|^2},\q k_{n,+,+}=\s{|x'-y'|^2+|x_3+y_3+2n|^2}.
\]

\noindent{\bf Case 1: $\boldsymbol{|x'-y'|<1}$}

\bea
G(x,y)&=&\f{1}{4\pi}\sum\limits^{+\i}_{n=-\i}\bigg\{\f{1}{k_{n,-,+}}-\f{1}{k_{n,+,-}}\bigg\}\nn\\
      &=&\f{1}{4\pi}\underbrace{\bigg\{\f{1}{k_{0,-,+}}-\f{1}{k_{0,+,-}}\bigg\}}_{I_1}+\f{1}{4\pi}\underbrace{\bigg\{\f{1}{k_{1,-,+}}-\f{1}{k_{1,+,-}}\bigg\}}_{I_2}\nn\\
      &&+\f{1}{4\pi}\sum\limits^{n\neq 0,1}_{n\in\bZ}\underbrace{\bigg\{\f{1}{k_{n,-,+}}-\f{1}{k_{n,+,-}}\bigg\}}_{I_{3,n}}.\nn
\eea
It is easy to see that
\bea\label{A.3}
|I_1|+|I_2|\ls  \f{1}{\s{|x'-y'|^2+|x_3-y_3|^2}}.
\eea
We compute $I_3$ as follows
\bea\label{A.4}
|I_{3,n}|&\ls&  \f{|k_{n,+,-}-k_{n,-,+}|}{k_{n,-,+}k_{n,+,-}}\nn\\
     &=&  \f{|k^2_{n,+,-}-k^2_{n,-,+}|}{k_{n,-,+}k_{n,+,-}(k_{n,-,+}+k_{n,+,-})}\nn\\
     &=&  \f{|4x_3(y_3-2n)|}{k_{n,-,+}k_{n,+,-}(k_{n,-,+}+k_{n,+,-})}.
\eea
When $n\in\bZ,n\neq 0,1$, we have
\be
|x_3+y_3-2n|\gtrsim |n|,\q  |x_3-y_3+2n|\gtrsim |x_3-y_3|+|n|, \nn
\ee
which indicate that
\be
k_{n,-,+}\gtrsim |x'-y'|+|x_3-y_3|+|n|,\q k_{n,+,-}\gtrsim |n|.\nn
\ee
Inserting the above inequality into \eqref{A.4}, we have
\bea
|I_{3,n}|&\ls&  \f{1}{|n|\s{|x'-y'|^2+|x_3-y_3|^2+n^2}}.\nn
\eea
So
\bea\label{A.5}
\sum\limits^{n\neq 0,1}_{n\in\bZ}|I_{3,n}|&\ls&\sum\limits^{n\neq 0,1}_{n\in\bZ}\f{1}{|n|\s{|x'-y'|^2+|x_3-y_3|^2+n^2}}\nn\\
&\ls&\int^\i_1 \f{1}{s\s{|x'-y'|^2+|x_3-y_3|^2+s^2}}ds\nn\\
&\ls& \f{1}{\s{|x'-y'|^2+|x_3-y_3|^2}}.
\eea
Then \eqref{A.3} and \eqref{A.5} together indicates the first part of \eqref{A.1}.

\noindent{\bf Case 2: $\boldsymbol{|x'-y'|>1}$}

Since in this situation, we need to get one more order decay with respect to $|x'-y'|$, we need to compute the sum more carefully.

\bea
G(x,y)&=&\f{1}{4\pi}\sum\limits^{+\i}_{n=-\i}\bigg\{\f{1}{k_{n,-,+}}-\f{1}{k_{n,+,-}}\bigg\}\nn\\
      &=&\f{1}{4\pi}\underbrace{\bigg\{\f{1}{k_{0,-,+}}-\f{1}{k_{0,+,-}}+\f{1}{k_{1,-,+}}-\f{1}{k_{1,+,-}}+\f{1}{k_{-1,-,+}}-\f{1}{k_{-1,+,-}}\bigg\}}_{J_1}\nn\\
      &&+\f{1}{4\pi}\sum\limits_{n\geq2}\underbrace{\bigg\{\f{1}{k_{n,-,+}}-\f{1}{k_{n,+,-}}
      +\f{1}{k_{n,-,-}}-\f{1}{k_{n,+,+}}\bigg\}}_{J_{2,n}}.\nn
\eea
Now we will compute $J$s term by term.
For $J_1$,
\bea\label{A.6}
\Big|\f{1}{k_{0,-,+}}-\f{1}{k_{0,+,-}}\Big|&=&\f{|k^2_{0,-,+}-k^2_{0,+,-}|}{{k_{0,-,+}}{k_{0,+,-}}({k_{0,-,+}}+{k_{0,+,-}})}\nn\\
                     &=&\f{|4x_3y_3|}{{k_{0,-,+}}{k_{0,+,-}}({k_{0,-,+}}+{k_{0,+,-}})}\nn\\
     &\ls&\f{1}{|x'-y'|^3}.
\eea
By the same techniques, we can prove that
\be\label{A.7}
\Big|\f{1}{k_{1,-,+}}-\f{1}{k_{1,+,-}}\Big|+\Big|\f{1}{k_{-1,-,+}}-\f{1}{k_{-1,+,-}}\Big|\ls\f{1}{|x'-y'|^3}.
\ee

We see that when $|x'-y'|>1$ and $n\geq 2$,
\be\label{A.8}
k_{n,\alp,\beta}\thickapprox |x'-y'|+n. \q \alp,\beta\in\{+,-\}.
\ee
Then
\bea
J_{2,n}&=&\f{k^2_{n,+,-}-k^2_{n,-,+}}{{k_{n,-,+}}{k_{n,+,-}}({k_{n,-,+}}+{k_{n,+,-}})}+\f{k^2_{n,+,+}-k^2_{n,-,-}}{{k_{n,-,-}}{k_{n,+,+}}({k_{n,-,-}}+{k_{n,+,+}})}\nn\\
     &=&\f{4x_3(y_3-2n)}{{k_{n,-,+}}{k_{n,+,-}}({k_{n,-,+}}+{k_{n,+,-}})}+\f{4x_3(y_3+2n)}{{k_{n,-,-}}{k_{n,+,+}}({k_{n,-,-}}+{k_{n,+,+}})}\nn\\
   &=&4x_3y_3\underbrace{\bigg[\f{1}{k_{n,-,+}k_{n,+,-}(k_{n,-,+}+k_{n,+,-})}+\f{1}{k_{n,-,-}k_{n,+,+}(k_{n,-,-}+k_{n,+,+})}\bigg]}_{K_{n,1}}\nn\\
   &&-8x_3n\underbrace{\bigg[\f{1}{k_{n,-,+}k_{n,+,-}(k_{n,-,+}+k_{n,+,-})}-\f{1}{k_{n,-,-}k_{n,+,+}(k_{n,-,-}+k_{n,+,+})}\bigg]}_{K_{n,2}}.\nn
\eea
So, we have
\bea
|J_{2,n}|\ls |K_{n,1}|+n|K_{n,2}|.\nn
\eea
Next we will show that
\bea\label{A.9}
\sum\limits^\i_{n=2} |J_{2,n}|&\ls& \sum\limits^\i_{n=2} |K_{n,1}|+ \sum\limits^\i_{n=2}n|K_{n,2}|\nn\\
                           &\ls&\f{1}{|x'-y'|^2}.
\eea
Then \eqref{A.6}, \eqref{A.7} and \eqref{A.9} together indicate the second part of \eqref{A.1}.

Using \eqref{A.8},
\bea
\sum\limits^\i_{n=2} |K_{n,1}|&\ls& \sum\limits^\i_{n=2} \f{1}{(|x'-y'|+n)^3}\nn\\
                              &\ls& \int^\i_0\f{1}{(|x'-y'|+s)^3}ds\ls  \f{1}{|x'-y'|^2}.\nn
\eea

The hardest part is to estimate $K_{n,2}$ since the sum has one more increasing term $n$ before $K_{n,2}$. When we estimate $K_{n,2}$, we need one more $\f{1}{n}$ coming out compared with $K_{n,1}$.  \\
Denote
\be
K_{n,2,1} \equiv {k_{n,-,+}\, \cdot k_{n,+,-} \cdot (k_{n,-,+}+k_{n,+,-})}\approx (|x'-y'|+n)^3,\nn
\ee
\be
K_{n,2,2} \equiv {k_{n,-,-} \, \cdot k_{n,+,+} \cdot (k_{n,-,-}+k_{n,+,+})}\approx (|x'-y'|+n)^3.\nn
\ee
So,
\bea\label{A.10}
|K_{n,2}|&=&\Big|\f{1}{K_{n,2,1}}-\f{1}{K_{n,2,2}}\Big|\nn\\
       &=&\f{|K^2_{n,2,2}-K^2_{n,2,1}|}{K_{n,2,1}K_{n,2,2}(K_{n,2,1}+K_{n,2,2})}\nn\\
       &\approx& \f{|K^2_{n,2,2}-K^2_{n,2,1}|}{(|x'-y'|+n)^9}.
\eea
By a direct computation, we can see that
\be\label{A.11}
|K^2_{n,2,2}-K^2_{n,2,1}|\ls |x'-y'|^4n+|x'-y'|^2n^3+n^5.
\ee
Inserting \eqref{A.11} into \eqref{A.10}, we can get
\bea
\sum\limits^\i_{n=2} n|K_{n,2}|&\ls& \sum\limits^\i_{n=2} \f{|x'-y'|^4n^2+|x'-y'|^2n^4+n^6}{(|x'-y'|+n)^9}\nn\\
          &\ls&\int^\i_0\f{|x'-y'|^4s^2+|x'-y'|^2s^4+s^6}{(|x'-y'|+s)^9}ds\nn\\
          &\ls& \f{1}{|x'-y'|^2}.
\eea

Combining the above, we have proved the estimate in \eqref{A.1} for $G(x,y)$. The estimate of $\p_{x',y'}G(x,y)$ will be essentially the same as $G(x,y)$, since one $\f{1}{|x-y|}$ will come out when we differentiate $G(x,y)$ on $x',y'$. So we omit the details.

{\bf Acknowledgement.} We wish to thank Professors Zhen Lei, Vladimir Sverak and Shangkun Weng, and Dr. Zijin Li for helpful conversations. X. H. Pan is supported by Natural Science Foundation of Jiangsu Province (No. SBK2018041027) and National Natural Science Foundation of China (No. 11801268).
Q. S. Zhang wishes to thank the Simons Foundation (grant 282153) for its support , and is grateful to Fudan University for its hospitality during his visit. Finally, thanks go to the anonymous referee for the careful evaluation and helpful corrections and suggestions.

\vskip 2cm

\qq \qq\qq\qq\qq\qq Bryan Carrillo,

\qq\qq\qq\qq\qq\qq Department of Mathematics,

\qq\qq\qq\qq\qq\qq University of California, Riverside, CA, 92521, USA;

\vskip 0.5cm

\qq \qq\qq\qq\qq\qq Xinghong Pan

\qq\qq\qq\qq\qq\qq Department of Mathematics,

\qq\qq\qq\qq\qq\qq Nanjing University of Aeronautics and Astronautics,

\qq\qq\qq\qq\qq\qq Nanjing 211106, China.

\vskip 0.5cm

\qq \qq\qq\qq\qq\qq Qi S. Zhang

\qq\qq\qq\qq\qq\qq Department of Mathematics,

\qq\qq\qq\qq\qq\qq University of California, Riverside, CA, 92521, USA;

\qq\qq\qq\qq\qq\qq \& School of Mathematical Sciences,

\qq\qq\qq\qq\qq\qq Fudan University, Shanghai, 200433, China.

\vskip 0.5cm

\qq \qq\qq\qq\qq\qq Na Zhao

\qq\qq\qq\qq\qq\qq  School of Mathematical Sciences,

\qq\qq\qq\qq\qq\qq Fudan University, Shanghai, 200433, China.

\qq\qq\qq\qq\qq\qq \& Department of Mathematics,

\qq\qq\qq\qq\qq\qq University of California, Riverside, CA, 92521, USA;

\end{document}